\documentclass[12pt]{amsart}
\font\emailfont=cmtt10

\usepackage{amsmath,amsthm,amsfonts,amscd,flafter,epsf}

\newcommand\commentable[1]{#1}

\newtheorem{theorem}{Theorem}[section]
\newtheorem{prop}[theorem]{Proposition}
\newtheorem{cor}[theorem]{Corollary}
\newtheorem{lemma}[theorem]{Lemma}

\newtheorem{defn}[theorem]{Definition}

\newtheorem{remark}[theorem]{Remark}

\def\endproof{\relax\ifmmode\expandafter\endproofmath\else
  \unskip\nobreak\hfil\penalty50\hskip.75em\hbox{}\nobreak\hfil\bull
  {\parfillskip=0pt \finalhyphendemerits=0 \bigbreak}\fi}
\def\endproofmath$${\eqno\bull$$\bigbreak}
\def\bull{\vbox{\hrule\hbox{\vrule\kern3pt\vbox{\kern6pt}\kern3pt\vrule}\hrule}}

\newcommand{\R}{\mathbb{R}}

\newcommand{\C}{\mathbb{C}}

\newcommand{\Z}{\mathbb{Z}}

\newcommand{\OneHalf}{\frac{1}{2}}

\newcommand{\ind}{\mathrm{ind}}

\newcommand{\cm}{\cdot}

\newcommand{\Nbd}[1]{\nu(#1)}
\newcommand{\nbd}[1]{\Nbd{#1}}
\newcommand{\CDisk}{D}

\newcommand{\Ideal}[1]{{\mathcal{I}}_{#1}}

\newcommand{\Sobol}[2]{L^{#1}_{#2}}

\newcommand{\Canon}[1]{{K}_{#1}}
\newcommand{\KSig}{\Canon{\Sigma}}

\newcommand{\Dirac}{\mbox{$\not\!\!D$}}

\newcommand{\Maps}{\mathrm{Map}}

\newcommand{\SpinC}{{\mathrm{Spin}}^c}
\newcommand{\Spin}{{\mathrm{Spin}}}
\newcommand{\tSpinC}{\mathrm{\underline{Spin}}^c}

\newcommand{\Proj}{\Pi}

\newcommand{\DDt}{\frac{\partial}{\partial t}}

\newcommand{\goesto}{\mapsto}

\newcommand{\DBar}{\overline{\partial}}

\newcommand{\SpecFlow}{\mathrm{SF}}
\newcommand{\SF}{\SpecFlow}

\newcommand\Wedge{\Lambda}

\newcommand\loc{\mathrm{loc}}

\newcommand\abuts\Rightarrow
\newcommand\Sym{\mathrm{Sym}}

\newcommand{\Hol}{\mathrm{Hol}}

\newcommand\cyl{cyl}
\newcommand\Cinfty{C^{\infty}}

\newcommand\Torus[1]{\mathbb{T}^{#1}}
\newcommand\Perm[1]{S_{#1}}
\newcommand\Met{{\mathfrak{Met}}}
\newcommand\Targ{\Torus{}(Y)}
\newcommand\TargP{{\mathbb T}'(Y)}
\newcommand\Lag{L}
\newcommand\Quot{Q}
\newcommand\Inv{\theta} 
\newcommand\OurJac{\widetilde{\Jac}}
\newcommand\OurSym{{\widetilde\Sym}^{g-1}(\Sigma)}
\newcommand\OurSymF{{\widetilde\Sym}^{g-1}_f(\Sigma)}
\newcommand\OurG{\Gamma}
\newcommand\OurTheta{\widetilde\Theta}

\newcommand\spinc{\mathfrak s}
\newcommand\spincX{\mathfrak r}

\newcommand\Jac{J}
\newcommand\spin{\mathfrak s}
\newcommand\NumEnds{n}
\newcommand\NumCurves{n}
\newcommand\Error{\epsilon}

\newcommand{\FFill}{F}

\newcommand{\EFill}{E}
\newcommand{\FCyl}{F^+}
\newcommand{\ECyl}{E^+}
\newcommand{\FComp}{F^c}
\newcommand{\Tors}{\mathrm{Tors}}

\newcommand{\met}{h} %{\mathfrak g}
\newcommand{\metk}{k} %{\mathfrak h}
\newcommand{\Conj}[1]{\overline{#1}}
\newcommand{\Wall}{\mathcal W}
\newcommand{\NdZero}{G}

\newcommand{\Splice}[1]{\gamma_{_{#1}}}
\newcommand{\Harm}{\mathcal H}

\newcommand{\ParTheta}{\Psi}

\newcommand{\AJ}{\mu}
\newcommand{\AbTa}[1]{{\widetilde{\mathbb T}}_{#1}(\alpha)}
\newcommand{\AbTb}[1]{{\widetilde{\mathbb T}}_{#1}(\beta)}

\newcommand{\AbAlpha}{\widetilde\alpha}
\newcommand{\AbBeta}{\widetilde\beta}
\newcommand{\Ta}[1]{{\mathbb T}_{#1}(\alpha)}
\newcommand{\Tb}[1]{{\mathbb T}_{#1}(\beta)}
\newcommand{\La}[1]{{\Lambda}_{#1}(\alpha)}
\newcommand{\Lb}[1]{{\Lambda}_{#1}(\beta)}
\newcommand{\BigO}{O}
\newcommand{\OurDelta}{\widetilde\delta}

\newcommand{\TorsFree}{\underline H}
\newcommand{\tH}{\underline H}
\newcommand{\tInv}{\underline\theta}

\newcommand{\MetCyl}{g_{\cyl}}
\newcommand{\MetComp}{g}
\newcommand{\DiracCyl}{\Dirac_{\cyl}}
\newcommand{\DiracDisk}{\Dirac}
\newcommand{\tDelta}{\underline\Delta}
\newcommand{\tC}{\underline C}

\title{The Theta Divisor and Three-Manifold Invariants}
\author[Peter Ozsv{\'a}th]{Peter Ozsv\'ath}
\thanks{The first author was partially supported by NSF grant number
9971950}
\address{Department of
Mathematics,  Princeton University, New Jersey 08540 \newline
\indent{\emailfont{petero@math.princeton.edu}}}

\author[Zolt{\'a}n Szab{\'o}]{Zolt{\'a}n Szab{\'o}} 
\thanks{The second author was partially 
supported by NSF grant number DMS 970435 and a
Sloan Fellowship} 
\address{Department of
Mathematics,  Princeton University, New Jersey 08540 \newline
\indent{\emailfont{szabo@math.princeton.edu}}}

\begin{document}

\begin{abstract}
In this paper we study an invariant for oriented three-manifolds with
$b_1>0$, which is defined using Heegaard splittings and the theta
divisor of a Riemann surface. The paper is divided into two parts, the
first of which gives the definition of the invariant, and the
second of which identifies it with more classical (torsion) invariants
of three-manifolds. Its close relationship with Seiberg-Witten theory
is also addressed.
\end{abstract}

\maketitle
\section{Introduction}
\label{sec:Introduction}

Let $Y$ be an oriented three-manifold whose first Betti number
$b_1(Y)>0$. In this paper, we study a topological invariant of $Y$,
which is a function $$\Inv\colon
\SpinC(Y)\longrightarrow \Z$$ 
on the set of $\SpinC$ structures on $Y$, defined using Heegaard
splittings. Roughly speaking, the invariant measures how the theta
divisor of a Riemann surface behaves under certain degenerations of
the metric which are naturally associated to the Heegaard splitting.
To facilitate a more precise description, we recall some relevant
objects associated to Riemann surfaces and then Heegaard splittings.

Fix a Riemannian surface $\Sigma$ of genus $g$. We think of the
Jacobian $\Jac$ as the space of complex line bundles ${\mathcal E}$
over $\Sigma$ of degree $g-1$, modulo isomorphism. A generic bundle
${\mathcal E}$ in $\Jac$, admits no holomorphic sections. The theta
divisor, then, is the locus of line bundles which do. Note that the
space $\Jac$ is a real $2g$-dimensional torus; indeed, a spin
structure naturally induces an identification between the space $\Jac$
and the torus $H^1(\Sigma;S^1)$. (Here, we think of the circle $S^1$
as $\R/\Z$.)  Moreover, the theta divisor is the image of the
Abel-Jacobi map $$\Theta\colon \Sym^{g-1}(\Sigma)
\longrightarrow \Jac,$$
which assigns to a divisor the corresponding holomorphic line bundle.

Now, consider a handlebody $U$ bounding $\Sigma$. Such a handlebody
gives rise to a canonical $g$-dimensional torus $L(U)$ in $\Jac$:
$L(U)$ corresponds to the image of $H^1(U;S^1)$ in $H^1(\Sigma;S^1)$
via the identification corresponding to a spin structure $\spinc_0$ on
$\Sigma$ which extends over $U$. (Note, however, that $L(U)$ is
independent of the choice of spin structure used in its definition.)

A handlebody $U$ bounding $\Sigma$ can be described using Kirby
calculus. $U$ is obtained from $\Sigma$ by first attaching $g$
two-handles along $g$ disjoint simple, closed curves
$\{\gamma_1,...,\gamma_g\}$ which are linearly independent in
$H_1(\Sigma;\Z)$; and then one three-handle. The collection
$\{\gamma_1,...,\gamma_g\}$ will be called {\em a complete set of
attaching circles} for $U$. Since the three-handle is unique, $U$ is
determined by a complete set of attaching circles.

A  handlebody $U$ gives rise to a class of $U$-allowable metrics
on $\Sigma$ (see Definitions~\ref{def:CurveAllowable} and
\ref{def:Allowable}), which correspond to certain degenerations of
$\Sigma$.  For instance, if $\{\gamma_1,...,\gamma_g\}$ is a complete
set of attaching circles for $U$, then any metric which is
sufficiently stretched out normal to all of the $\gamma_i$ is
$U$-allowable. One special property of a $U$-allowable metric is
that the corresponding theta divisor is always disjoint from the
subspace $L(U)$ (see Lemma~\ref{lemma:MissTheta}).

Recall that a genus $g$ Heegaard decomposition of an oriented $3$-manifold 
is a decomposition of
$Y=U_0\cup_{\Sigma} U_1$ into two handlebodies $U_0$ and $U_1$ which
are identified along their boundary, which is a surface $\Sigma$ of
genus $g$. Denote by $L_i$ the associated tori $L(U_i)$ in the
Jacobian.  Fix a one-parameter family $h_t$ of metrics on $\Sigma$ for
which $h_0$ is $U_0$-allowable, and $h_1$ is $U_1$-allowable. Then,
consider the set of points in
$[0,1]\times[0,1]\times\Sym^{g-1}(\Sigma)$ $$\{(s,t,D)\big|
s\leq t \text{~and~}
\Theta_{h_s}(D)\in L_0 \text{~and~}
\Theta_{h_t}(D)\in L_1\}.$$
We show that for small, generic perturbations of $L_i$, this set of
points is isolated. Moreover, there is a natural map from this set to
the set of $\SpinC$ structures on $Y$. Then, $\theta(\spinc)$ is a
signed count of the number of points corresponding to the $\SpinC$
structure $\spinc$.

A geometric meaning of this signed count can be given as follows. The
one-parameter family of metrics induces a map $$\Theta \colon
[0,1]\times \Sym^{g-1}(\Sigma)\rightarrow \Jac,$$ by
$\Theta(t,D)=\Theta_{\met_t}(D)$. The set $\Theta^{-1}(L_0)$ misses
the region where $t<\epsilon$, and $\Theta^{-1}(L_1)$ misses the
region where $t>1-\epsilon$. The tori $L_0$ and $L_1$ can be perturbed
slightly to make them disjoint.  The invariant $\theta$ then measures
the degree to which the preimages under $\Theta$ of these perturbed
versions of $L_0$ and $L_1$ are linked. One gets more than a simple
linking number -- hence the function on $\SpinC(Y)$ -- by passing to a
suitable covering space of $\Jac$, and looking at the linking numbers
between the preimages of the various lifts of $L_0$ and $L_1$. Details
are spelled out in Section~\ref{sec:DefTheta}, where the first main
result is the following:

\begin{theorem}
\label{thm:WellDefined}
The invariant
$$\Inv\colon
\SpinC(Y)\longrightarrow \Z$$ 
is well-defined; in particular, it does not depend on the metrics,
perturbations, and Heegaard decompositions of $Y$.
\end{theorem}

The invariant $\Inv$ also manifestly shares some of the properties of the
Seiberg-Witten invariant for three-manifolds.

\begin{prop}
\label{prop:StandardProperties}
For any given oriented three-manifold $Y$, with $b_1(Y)>0$ 
there are only finitely many
$\SpinC$ structures for which $\theta(\spinc)\neq 0$.  Moreover,
$\theta(\spinc)=\theta(\Conj{\spinc})$, where the map $\spinc\mapsto
\Conj{\spinc}$ denotes the natural involution on the set
$\SpinC(Y)$. Also, if $-Y$ denotes the oriented manifold obtained by
reversing the orientation of $Y$, then
$$\theta_Y(\spinc)=(-1)^{b_1+1}\theta_{-Y}(\spinc).$$
\end{prop}

After laying down the basis for the definition of the invariant, we
turn to its computation. It will be convenient for us to think of the
invariant as an element $\theta\in\Z[\SpinC(Y)]$, in the usual manner:
$$\Inv=\sum_{\spinc\in\SpinC(Y)}\Inv(\spinc)[\spinc],$$
where $\Z[\SpinC(Y)]$ is to be thought of as a module over the group-ring
$\Z[H]$ associated to the group $H=H^2(Y;\Z)\cong H_1(Y;\Z)$.
In fact, in the computation, we begin by considering a
weaker invariant, obtained from $\theta$ by dividing out by the action
of the torsion subgroup $\Tors$ of $H^2(Y;\Z)$. In keeping with the
convention of \cite{MengTaubes}, we underline objects when they are to
be viewed modulo the action of the torsion subgroup $\Tors$ so,
e.g. $\tH$ and $\tSpinC(Y)$ denote the quotients
of $H$ and $\SpinC(Y)$ respectively by the
action of $\Tors$. There is 
an induced invariant 
$$\tInv\colon \tSpinC(Y)
\longrightarrow \Z,$$ 
defined by adding the values of $\theta(\spinc)$ for all $\SpinC$
structures in a given orbit.
Via the natural
identification 
\begin{equation}
\label{eq:IdentTSpinC}
\tSpinC(Y)\cong\tH,
\end{equation}
which sends any $\Spin$ structure to $0$,
we can view $\tInv$ as an element $\tInv\in\Z[\tH]$.

\begin{theorem}
\label{thm:Calculate}
If $b_1(Y)>1$, then up to sign, ${\underline \theta}$ is equal to the
symmetrized Alexander polynomial of $Y$.
\end{theorem}

\begin{theorem}
\label{thm:CalculateZ}
Suppose $b_1(Y)=1$, and let $A=a_0+\sum_{i=1}^k a_i(t^i+t^{-i})$ be
the symmetrized Alexander polynomial of $Y$ normalized so that
$A(1)=|\Tors H_1(Y;\Z)|$. Then,
$${\underline\theta}(i)=\sum_{j=1}^\infty j\cdot a_{|i|+j},$$ (note
that we are using the natural identification $\tSpinC(Y)\cong \Z$
coming from ~\eqref{eq:IdentTSpinC}).
\end{theorem}

In fact, a closer inspection of the proofs of
Theorems~\ref{thm:Calculate} and \ref{thm:CalculateZ} gives a more
refined statement, which identifies the invariant $\theta$ with a
torsion invariant $\tau\in\Z[\SpinC(Y)]$ discovered by Turaev,
see~\cite{Turaev}. (The element we denote by $\tau$ here is the element of
$\Z[\SpinC(Y)]$ induced from Turaev's ``torsion function'' $T$ of \S5
from~\cite{Turaev}.)

\begin{theorem}
\label{thm:CalcBigT}
Suppose $b_1(Y)>1$. Then the invariant $\theta\in\Z[\SpinC(Y)]$
agrees, up to possibly translation by two-torsion in $H^2(Y;\Z)$ and a
sign which depends only on $b_1(Y)$, with the Turaev invariant $\tau$.
\end{theorem}

When $b_1(Y)=1$, Turaev's torsion function depends on a choice of
generator of $\tH$. We recall from Turaev that if one fixes an
$\tH\cong \Z$ and $t$ denotes the positive generator of $\tH$, then
the two torsion functions $T_t$ and $T_{t^{-1}}$ are related by the
formula
$$T_{t^{-1}}(\spinc)=T_t(\spinc) - {\underline\spinc}.$$
Moreover, the support of $T_t$ and $T_{t^{-1}}$ are bounded above and
below respectively. Now, one can define a compactly supported torsion
function $T'$ which does not depend on a choice of generator by
$T'(\spinc)=T_t(\spinc)$ if ${\underline\spinc}$ is a non-negative
multiple of $t$, $T'(\spinc)=T_{t^{-1}}(\spinc)$ otherwise; or
equivalently
$$T'(\spinc)
=\frac{1}{2}(T_t(\spinc)+T_{t^{-1}}(\spinc)+|{\underline\spinc}|).$$
Let $$\tau'=\sum_{\spinc}T'(\spinc)[\spinc].$$
Our result can then be stated as follows:

\begin{theorem}
\label{thm:CalcOneT}
Suppose $b_1(Y)=1$. Then the invariant $\theta\in\Z[\SpinC(Y)]$ agrees,
up to possibly translation by two-torsion in $H^2(Y;\Z)$, with the
Turaev invariant $\tau'$.
\end{theorem}

The relationship between the invariant $\Inv$ and the Seiberg-Witten
invariant for three-manifolds can be seen from two different points of
view. On the one hand, the invariant arises naturally when studying
the Seiberg-Witten equations for Heegaard decompositions; in fact this
is how we discovered it. On the other hand, results of
Meng-Taubes~\cite{MengTaubes} and Turaev~\cite{TuraevTwo}, together
with our computation, show that the invariant $\Inv$ agrees with a
numerical invariant obtained from the Seiberg-Witten equations. It is
also interesting to compare this with the Morse-theoretic
constructions of~\cite{HutchLee}, and also with recent work of
Salamon~\cite{SalamonSW}.

Throughout the paper, we work with three-manifolds whose first Betti
number is positive. In the case where $b_1(Y)=0$, there is a naturally
associated invariant, which is technically more complicated to
describe. The reason for this is that, when $b_1(Y)=0$, the invariant
$\theta$ actually depends on the path of metrics used in its
definition, so to get a topological quantity, one must correct by a
spectral flow correction term. These issues are addressed
in~\cite{ThetaCasson}, where the relationship between this construction
and the Casson-Walker invariant is explored.

The present paper is organized as follows. Roughly speaking, it can be
divided into two parts: the first of which
(Sections~\ref{sec:DefTheta}--\ref{sec:Splicing}, together with
Section~\ref{sec:WallCross}) defines the invariant, and the second of
which (Sections~\ref{sec:CalcBig}--\ref{sec:Alex}) calculates it.  In
Section~\ref{sec:DefTheta}, we describe the metrics on $\Sigma$
induced by the Heegaard decomposition of $Y$, and show that the
invariant $\theta(\spinc)$ is independent of choices of metrics, and
hence can depend only on the Heegaard decomposition of $Y$ (except for
the special case where $Y$ is a rational homology $S^1\times S^2$ but
not an an integer homology $S^1\times S^2$, a case which we return to
in Section~\ref{sec:WallCross}). The results rely on a few technical
lemmas about the behaviour of the theta divisor under degenerations of
$\Sigma$, which are proved in
Section~\ref{sec:MoveTheta}. Independence of the Heegaard
decomposition, then, amounts to proving ``stabilization invariance''
of the invariant. This result is proved in Section~\ref{sec:Splicing},
as a corollary to some results about the behaviour of the theta
divisor under degenerations of the metric along homologically
inessential curves.  The degenerations of Section~\ref{sec:Splicing}
play an important role in the second part of the paper, as well.  The
calculation of the invariant depends on a certain perturbation, which
involves slightly enlarging the tori $L(U)\subset
\Jac$ coming from the handlebodies, to a $g+1$-dimensional torus which
intersects the theta divisor even when the metric on $\Sigma$ is
$U$-allowable. Another corollary of the results of
Section~\ref{sec:Splicing}, then, is an explicit understanding of the
intersection of the theta divisor with these larger tori.

With the technical background in place, we turn to the
calculations, which identify the invariant $\theta$ with data of a more
directly topological character. In Section~\ref{sec:CalcBig}, we focus
on the case when $b_1(Y)>1$, which is slightly simpler than
the calculation when
$b_1(Y)=1$ given in Section~\ref{sec:CalcOne}, as there is
more freedom in perturbing the invariant when the second Betti number
is large. But the same general idea works in both cases. The close
relationship between the topological data obtained and 
the Alexander polynomial (see Theorems~\ref{thm:Calculate} and
\ref{thm:CalculateZ} above) is explained in
Section~\ref{sec:Alex}. Indeed, a closer look at the proofs of these
results gives the more refined formulations involving Turaev's torsion
invariant (see Theorems~\ref{thm:CalcBigT} and \ref{thm:CalcOneT}
above), as shown in Subsection~\ref{subsec:Refinements}.

A final debt is paid in Section~\ref{sec:WallCross}, where we address
the case of topological invariance in the case where $b_1(Y)=1$ (with
no assumptions on the torsion in $H_1(Y;\Z)$). This section should be
thought of as an appendix to the first part of the paper, though the
proofs are of a slightly different character than those in the rest of
the paper (bearing a closer relationship to the issues addressed
in~\cite{ThetaCasson}).

\section{Defining the invariant}
\label{sec:DefTheta}

The aim of this section is to spell out the details that go into the
definition of the invariant $\theta$ sketched in
Section~\ref{sec:Introduction}, and, indeed to prove that its value
depends only on the topology of the Heegaard decomposition of the
three-manifold. We complete the proof of Theorem~\ref{thm:WellDefined}
in Section~\ref{sec:Splicing}, where we prove, among other things,
that $\theta$ remains invariant under stabilization.

Fix an oriented three-manifold $Y$ whose first Betti number
$b_1(Y)>0$. The definition of the invariant makes reference to a
genus $g$ Heegaard decomposition of $Y=U_0\cup_{\Sigma}U_1$; so we
discuss some objects naturally associated to such a decomposition.

We give first a convenient definition of the Jacobian of an oriented
$2$-manifold $\Sigma$ of genus $g$, endowed with a Riemannian metric
$\met$. Fix a Hermitian line bundle $E$ over $\Sigma$ whose Euler
number is $g-1$. A Hermitian connection $A$ over $E$ is said to have
{\em normalized curvature form} if its curvature form satisfies $$F_A
=
\OneHalf F_{K(\met)},$$ where $K(\met)$ is the Levi-Civita connection on the
canonical bundle for the metric $\met$.  Then the Jacobian $\Jac_\met$
is the space of Hermitian connections $A$ with normalized curvature
form, modulo the gauge group of circle-valued functions
$\Maps(\Sigma,S^1)$.  The group $H^1(\Sigma;\R)$ acts simply
transitively on $\Jac_\met$, with stabilizer $H^1(\Sigma;\Z)$ and so a
point in $\Jac_\met$ gives an identification of $\Jac_\met$ with the
$2g$-dimensional torus $$\Jac_\met\cong
\frac{H^1(\Sigma;\R)}{H^1(\Sigma;\Z)}=H^1(\Sigma;S^1).$$
Moreover, a spin structure on $\Sigma$ naturally gives rise to a point
in $\Jac_\met$, and hence an identification $$\Jac_\met\cong
H^1(\Sigma;S^1).$$ When it is clear from the context, we drop the
metric $\met$ from the notation for the Jacobian.

We will typically work in a certain cover of the Jacobian which is
associated to the Heegaard decomposition. Specifically, the long exact
sequence in cohomology for the decomposition induces a (surjective)
coboundary map $\delta\colon H^1(\Sigma;\Z)\longrightarrow H^2(Y;\Z)$,
whose kernel we denote by $\OurG$ (alternatively, this is the subgroup
of $H^1(\Sigma;\Z)$ generated by the image of $H^1(U_0;\Z)\oplus
H^1(U_1;\Z)$ under the obvious inclusion map).  We find it convenient,
then, to consider the cover $\OurJac$ of $\Jac$, the space of
connections in $E$ with normalized curvature form, modulo the action
by gauge transformations in the kernel of the composite $$\begin{CD}
\Maps(\Sigma;S^1)@>>> H^1(\Sigma;\Z) @>{\delta}>> H^2(Y;\Z)
\end{CD}
$$ The space $\OurJac$ inherits a natural action of $H^1(\Sigma;\R)$,
and the action of $H^1(\Sigma;\Z)$ on $\OurJac$ descends to a free
action of $$\frac{H^1(\Sigma;\Z)}{\OurG}\cong H^2(Y;\Z)$$ on
$\OurJac$, whose quotient is canonically identified with $\Jac$.  The
condition that $b_1(Y)>0$ is equivalent to the condition that
$\OurJac$ is a non-compact space.

Given a metric $\met$, there is an ``Abel-Jacobi map''
$$\Theta_\met\colon \Sym^{g-1}(\Sigma)\longrightarrow
\Jac_\met,$$
where $\Sym^{g-1}(\Sigma)$ is the space of effective degree $g-1$
divisors on $\Sigma$, i.e. $g-1$-fold symmetric power of $\Sigma$
(note that our conventions are slightly different from those typical
in Riemann surface theory, where the Jacobian is often thought of as
the group of complex structures on a topologically trivial line
bundle, rather than the positive spinor bundle).  Given a divisor
$D\in
\Sym^{g-1}(\Sigma)$, the corresponding connection $\Theta_\met(D)$ is
characterized by its curvature form (half that of the canonical bundle
with metric $\met$) and its associated $\DBar$-operator, which we
require to admit a holomorphic section which vanishes exactly at
$D$. The image of this map in $\Jac$ is called the {\em theta
divisor}. Once, again, we find it convenient to work in a lift
$\OurSym$ of $\Sym^{g-1}(\Sigma)$. This lift corresponds to the
subgroup of $\pi_1(\Sym^{g-1}(\Sigma))$ which is the kernel of the
composite $$\begin{CD}
\pi_1(\Sym^{g-1}(\Sigma))@>>>
H_1(\Sym^{g-1}(\Sigma))@>{(\Theta_\met)_*}>>
H_1(\Jac)\cong H^1(\Sigma;\Z) @>{\delta}>> H^2(Y;\Z),
\end{CD}
$$ where the first map is the Hurewicz homomorphism. Clearly,
$(\Theta_{\met})_*$ is independent of the Riemannian metric.  Thus, we
have a map $$\OurTheta_\met\colon \OurSym\longrightarrow \OurJac$$
which fits into a commutative diagram $$\begin{CD}
\OurSym @>{\OurTheta_\met}>> \OurJac \\
@VVV  @VVV \\
\Sym^{g-1}(\Sigma) @>{\Theta_\met}>> \Jac.
\end{CD}
$$

Note that standard Hodge theory gives an identification between the
Jacobian and the theta divisor given here with the definitions used in
the introduction.

Fix a handlebody $U$ which bounds $\Sigma$, and view the group
$H^1(U;\R)$ as a subgroup of $H^1(\Sigma;\R)$ using the natural
inclusion.  There is a natural quotient map
$$\Quot_U\colon\Jac\longrightarrow
\frac{H^1(\Sigma;S^1)}{H^1(U;S^1)}\cong H^2(U,\Sigma;S^1),$$ 
given as follows. Fix
a $\Spin$ structure $\spin_0$ on $U$, and let $p\in\Jac$ be the
induced point in the Jacobian. Given any $B\in\Jac$, there is a unique
$a\in H^1(\Sigma;S^1)$ so that $B=p+a$; we define $\Quot_U(B)$ to be
$a$ (modulo $H^1(U;S^1)$). This coset is independent of the spin
structure on $U$ since any two spin structures on $U$ differ by a
translation by a cohomology class coming from $H^1(U;S^1)$. The torus
$\Lag(U)$ defined in the introduction, then, is the preimage
$\Quot_U^{-1}(0)$.  Given a point $B\in\Jac$ and a homology class
$[\gamma]\in H_1(\Sigma;\Z)$ which bounds in $U$ there is a
well-defined holonomy, $\Hol_{\gamma}(B)\in S^1$, which is the
Kronecker pairing of $\Quot_U(B)$ with $[\gamma]$.

For a Heegaard decomposition of $Y$, let $\Lag_0$ and $\Lag_1$ denote
the associated tori $\Lag(U_0)$ and $\Lag(U_1)$ in $\Jac$.  A $\SpinC$
structure $\spinc$ on $Y$ gives rise to a pair of $g$-dimensional tori
$L_0(\spinc)$ and $L_1(\spinc)$ in $\OurJac$, up to simultaneous
translation by $H^2(Y;\Z)$, as follows. Let $\spin_0$ be a spin
structure on $Y$ and let $p$, $L_0$ and $L_1$ be as above.  Any
$\SpinC$ structure $\spinc$ on $Y$ can be written as $\spin_0+\ell$,
where $\ell\in H^2(Y;\Z)$.  Let ${\widetilde p}$ be any lift of $p$ to
$\OurJac$. Then, $L_0(\spinc)$ is the lift of $\Lag_0$ to $\OurJac$
which passes through ${\widetilde p}$, and $L_1(\spinc)$ is the lift
of $\Lag_1$ which passes through ${\widetilde p}+\ell$ (the translate
of ${\widetilde p}$ by the natural action of $H^2(Y;\Z)$ on
$\OurJac$).  Once again, it is easy to see that the subspaces are
independent of the spin structure $\spin_0$.  Since the intersection
in $H^1(\Sigma;\R)$ of the image of $H^1(U_0;\R)$ with $\OurG$ is
$H^1(U_0;\Z)$, it follows that $L_0(\spinc)$, and similarly
$L_1(\spinc)$, are both $g$-dimensional tori embedded in $\OurJac$.  A
{\em torsion $\SpinC$ structure} is a $\SpinC$ structure whose
associated real cohomology class ${\underline\spinc}=0$ (this is
equivalent to the condition that its first Chern class is torsion).
Note that $\spinc$ is torsion if and only if $L_0(\spinc)$ and
$L_1(\spinc)$ intersect. In fact, if $\spinc$ is torsion, then
$L_0(\spinc)\cap L_1(\spinc)$ is identified with $H^1(Y;\R)/H^1(Y;\Z)$.

Having introduced the basic topological objects associated to a Heegaard
decomposition, we must flesh out the notion of allowable metrics
used in the definition of the invariant $\theta$. The definition
corresponds to degenerations of the metric on $\Sigma$. We describe
these presently.

Let $\{\gamma_1,...,\gamma_{\NumCurves}\}$ be a collection of disjoint
simple, closed curves in $\Sigma$. Choose a tubular neighborhood $\nu$
of $\coprod_{i=1}^g\gamma_i$, and let $\met$ be a metric which extends
the product metric over $\nu$ arising naturally from an identification
$$\nu\cong
\coprod_{i=1}^g [-1,1]\times S^1.$$ Such a metric will be called {\em
product-like near the $\gamma_i$}. Given such a metric, let
$\met(T_1,...,T_{\NumCurves})$ denote the metric
obtained by inserting a tube of length $2T_i$ around the curve
$\gamma_i$, i.e. $\met(T_1,...,T_\NumCurves)$ is obtained by attaching
$$\coprod_{i=1}^g [-T_i,T_i]\times S^1$$ to $(\Sigma-\nu,\met)$ in the
obvious manner. The following lemma, whose proof is given in
Section~\ref{sec:MoveTheta}, describes what happens to the theta
divisor as the metric is stretched normal to curves in this way.

\begin{lemma} 
\label{lemma:MissTheta}
Let $U$ be a handlebody with boundary $\Sigma$, and let
$\{\gamma_1,...,\gamma_g\}$ be a complete set of attaching circles for
$U$. Then, for any compact set of metrics ${\mathcal H}$ on $\Sigma$
which are product-like near the $\gamma_i$, as the metrics are
stretched out normal to the $\gamma_i$, the theta divisor converges as
a point set, into $$\Hol_{\gamma_1}^{-1}(\OneHalf)\cup ... \cup
\Hol_{\gamma_g}^{-1}(\OneHalf);$$
i.e. given any $\epsilon$, there is a $T_0$ so that for all metrics
$h\in{\mathcal H}$, and for all $g$-tuples $(T_1,...,T_g)$ for which each
$T_i>T_0$, we have that $$
\Theta_{\met(T_1,...,T_g)}(\Sym^{g-1}(\Sigma))
\subset \bigcup_{i=1}^g \Hol_{\gamma_i}^{-1}
(\OneHalf-\epsilon, \OneHalf + \epsilon).$$
\end{lemma}

Note that $L(U)$ could be described as the set of points
$B\in\Jac(\Sigma)$ with $\Hol_{\gamma_i}B=0$ for all $i$. 
Thus, the above lemma says that for all metrics which are sufficiently
stretched out normal to all the $\gamma_i$, the theta divisor misses the
torus $L(U)$.  Indeed, it allows us to identify a special
(path-connected) class of metrics.

\begin{defn}
\label{def:CurveAllowable}
Let $U$ be a handlebody, and let $\{\gamma_1,...,\gamma_g\}$ be a complete
set of attaching circles for $U$. Fix a metric $\metk_0$
on $\Sigma$ which is sufficiently stretched out normal to the
$\{\gamma_i\}$ according to Lemma~\ref{lemma:MissTheta}. Another
metric $\metk_1$ on $\Sigma$ is called {\em allowable for
$\{\gamma_1,...,\gamma_g\}$}, or simply
$\{\gamma_1,...,\gamma_g\}$-allowable if there is a path $\metk_t$
connecting $\metk_0$ to $\metk_1$, so that
$$\Theta_{\metk_t}(\Sym^{g-1}(\Sigma))\cap L(U)=\emptyset$$
\end{defn}

In fact, Lemma~\ref{lemma:MissTheta} shows that the notion of
allowable is independent of the fixed metric $\metk_0$.  It appears,
however, to depend on the choice of the
$\{\gamma_1,...,\gamma_g\}$. The following proposition shows that this
is not the case: the notion of allowable depends only on the
handlebody $U$:

\begin{prop}
\label{prop:AllowableIsOK}
The class of allowable metrics depends only on the handlebody $U$;
i.e. if $\{\gamma_1,...,\gamma_g\}$ and $\{\gamma_1',...,\gamma_g'\}$
are any two complete sets of attaching circles for $U$, then a metric
is $\{\gamma_1,...,\gamma_g\}$-allowable if and only if it is
$\{\gamma_1',...,\gamma_g'\}$-allowable.
\end{prop}

The proof relies on the following lemma, whose proof is given in
Section~\ref{sec:MoveTheta}.

\begin{lemma}
\label{lemma:MissThetaHS}
Let $U$ be a handlebody with boundary $\Sigma$, and  let
$\{\gamma_1,...,\gamma_g\}$ be a complete set of attaching circles for
$U$. Then, for any compact set of metrics
${\mathcal H}$ on $\Sigma$ which are product-like near the $\gamma_i$
for $i=1,...,g-1$, given any $\epsilon>0$, there is a $T_0$ so that
for all $g-1$ tuples $(T_1,...,T_{g-1})$ with $T_i\geq T_0$ for all
$i$, and all metrics $\met\in{\mathcal H}$, we have that
$$\Theta^{-1}_{\met(T_1,...,T_{g-1})}(\Hol_{\gamma_1}^{-1}(0)\cap...\cap
\Hol_{\gamma_{g-1}}^{-1}(0))\subset \Hol_{\gamma_g}^{-1}(\OneHalf-\epsilon,\OneHalf+\epsilon).$$
\end{lemma}

The above lemma implies the following corollary:

\begin{cor}
\label{cor:OneCurveAtATime}
Let $\Sigma$, $\{\gamma_1,...,\gamma_g\}$, and $U$ be as in
Lemma~\ref{lemma:MissTheta}. For any  metric $\met$ which is
product-like around $\{\gamma_1,...,\gamma_{g-1}\}$, there is a constant
$T_0$ so that for all collections $T_i\geq T_0$, the metric
$\met(T_1,...,T_{g-1})$ is $\{\gamma_1,...,\gamma_g\}$-allowable. 
\end{cor}

\begin{proof}
Fix an initial metric $\metk$ which is product-like along all
$\{\gamma_1,...,\gamma_g\}$, and which agrees with $\met$ away from a
tubular neighborhood of $\gamma_g$. Lemma~\ref{lemma:MissTheta} gives
us a constant $C_0$ so that for all $g$-tuples $(T_1,...,T_g)$ with
$T_i\geq C_0$, $\metk(T_1,...,T_g)$ is allowable.  We can view the
metric $$\met_t(T_1,...,T_{g-1}) = t
\met(T_1,...,T_{g-1}) + (1-t)
\metk(T_1,...,T_{g-1},C_0)$$ as the result of inserting 
tubes with
parameters $T_1,...,T_{g-1}$ into a one-parameter (compact) family of
metrics away from the $\{\gamma_1,...,\gamma_{g-1}\}$. Thus,
Lemma~\ref{lemma:MissThetaHS} gives us a number $C_1$ with the
property that if all $T_i\geq C_1$, then for all metrics in the path
$\met_t(T_1,...,T_{g-1})$, the theta divisor misses $L(U)$. Hence,
if $T_1,...,T_{g-1} \geq\max(C_0,C_1)$, then
$\met(T_1,...,T_{g-1})$ is a $\{\gamma_1,...,\gamma_g\}$-allowable
metric.
\end{proof}

Proposition~\ref{prop:AllowableIsOK}, then, follows easily:

\vskip0.3cm
\noindent{\bf{Proof of Proposition~\ref{prop:AllowableIsOK}.}}
Fix $U$, and let $\{\gamma_1,...,\gamma_g\}$,
$\{\gamma_1',...,\gamma_g'\}$ be two complete sets of attaching
circles.  By standard Kirby calculus~\cite{KirbyCalc}, we see that it is
always possible to move between any two collections
$\{\gamma_1,...,\gamma_g\}$ and $\{\gamma_1',...,\gamma_g'\}$, through
a sequence of handle-slides.  Since a handle-slide fixes $g-1$ of the
curves, Corollary~\ref{cor:OneCurveAtATime} shows that the notion of
allowable remains unchanged.
\qed
\vspace{0.2in}

Proposition~\ref{prop:AllowableIsOK} allows us to refine the earlier
definition of allowable metrics:

\begin{defn}
\label{def:Allowable}
Let $U$ be a handlebody. A metric $\metk$ on $U$ is called {\em
$U$-allowable} provided that there is a complete set of attaching circles
$\{\gamma_1,...,\gamma_g\}$ for which the metric is allowable.
\end{defn}

With this background in place, we now give a definition of
$\theta(\spinc)$, where $\spinc\in\SpinC(Y)$.
Fix a smooth path of metrics $\{\met_t\}_{t\in[0,1]}$ for which $\met_0$ is $U_0$-allowable,
and $\met_1$ is $U_1$-allowable.
Consider the smooth map
$$\ParTheta\colon \Sym^{g-1}(\Sigma)\times\{(s,t)\in [0,1]\times[0,1] \big|
s\leq t\}\longrightarrow \frac{H^1(\Sigma;S^1)}{H^1(U_0;S^1)}
\times
\frac{H^1(\Sigma;S^1)}{H^1(U_1;S^1)}=\Targ,$$ defined by
$$\ParTheta(D,s,t)=\Quot_0(\Theta_{\met_s}(D))\times
\Quot_1(\Theta_{\met_t}(D)),$$ where $\Quot_i$ denotes the quotient map $\Quot_{U_i}$ for $i=0,1$, 
and let $M_{\eta_0\times \eta_1}$ denote the pre-image under
$\ParTheta$ of the point $\eta_0\times \eta_1 \in
\Targ$.
By Sard's theorem, for generic $\eta_0\times\eta_1$, this fiber is a
compact, oriented zero-dimensional manifold which misses the locus of
points $(D,s,t)$ where $s=t$.  Moreover, the points in
$M_{\eta_0\times\eta_1}$ can naturally be partitioned into subsets
indexed by the various $\SpinC$ structures on $Y$. Specifically, for a
choice of $\SpinC$ structure on $Y$ and corresponding lifts
$\Lag_i(\spinc)$ of $\Lag_i$ (for $i=0,1$), there is a subset of
$M_{\eta_0\times\eta_1}$, denoted $M_{\eta_0\times\eta_1}(\spinc)$,
corresponding to points $$\{(D,s,t)\in\OurSym\times[0,1]\times[0,1]
\big| \OurTheta_{\met_s}(D)\in \Lag_0(\spinc)+\eta_0, \OurTheta_{\met_t}(D)\in
\Lag_1(\spinc)+\eta_1, s\leq t\}$$ (where, in the above expression, $L_i(\spinc)+\eta_i$ denotes
the translate of the subset $L_i(\spinc)$ by $\eta_i$). Then,
$\theta_{\eta_0\times\eta_1}(\spinc)$ is defined to be the signed
number of points in this subset.  In the next two propositions, we
shall see that (given $\met_t$) there is an open neighborhood
$\NdZero$ of zero in $\Targ$ with the property that
$\theta_{\eta_0\times\eta_1}$ is independent of the particular
(generic) choice of $\eta_0\times\eta_1\in \NdZero$. This is
technically somewhat easier when $b_1(Y)>1$, so we consider that case
first. But before we do that, we pause for a moment to discuss signs.

Since $\Sym^{g-1}(\Sigma)\times[0,1]\times[0,1]$ is naturally
oriented, the sign of $\theta$ is determined by an orientation for the
torus $$\frac{H^1(\Sigma;S^1)}{H^1(U_0;S^1)} \times
\frac{H^1(\Sigma;S^1)}{H^1(U_1;S^1)},$$ which in turn is determined by
an ordering of the attaching circles $\{\alpha_1,...,\alpha_g\}$ and
$\{\beta_1,...,\beta_g\}$. We use an ordering for these attaching
circles which is consistent with the orientation of $H_*(Y)$, arising
from Poincar\'e duality, in the following sense.  We explain what this consistence means. The
Heegaard decomposition gives a chain complex for $Y$ with one
zero-cell, $g$ one-cells in one-to-one correspondence with the circles
$\{\alpha_i\}$, $g$ two-cells which correspond to the $\{\beta_i\}$,
and one three-cell.  In general, an orientation for a (finite
dimensional) chain complex is canonically equivalent to an orientation
for its real homology, since there is a splitting: $$\bigoplus_i C_i =
\bigoplus_i \left((\partial B_{i+1})\oplus H_i
\oplus B_i\right),$$ where $B_i\subset C_i$ is a vector space which is
mapped isomorphically under the boundary homomorphism to the group of
boundaries in $C_{i-1}$, which givesa natural identification of
$$\bigoplus_i C_i = \left(\bigoplus_i H_i\right)
\oplus \left(\bigoplus_i
B_i \oplus (\partial B_i) \right)$$ as oriented vector spaces; and the
vector space $\bigoplus_i (B_i \oplus (\partial B_i))$ is canonically
oriented as follows: if $\{b_j\}$ is a basis for $\bigoplus_i B_i$, we
declare that $$\bigoplus_j (b_j \oplus \partial b_j)$$ is a positive
oriented basis for $\bigoplus_i (B_i \oplus (\partial B_i))$.  The
orientation on $Y$ gives a canonical orientation of $H_0=C_0$,
$H_3=C_3$. Thus, in light of the above remarks, the orientation of
$H_*(Y)$ induces an orientation of $C_1\oplus C_2$, and hence an
ordering of the attaching circles.

Having nailed down the sign, we state the result we have been aiming
for, first in the case where $b_1(Y)>1$. 

\begin{prop}
\label{prop:IndepOfMetBig}
When $b_1(Y)>1$, the invariant $\theta(\spinc)$ depends only on the
Heegaard decomposition of $Y$; i.e. it is independent of the metrics
and perturbations used. More precisely, given any one-parameter family
$\met_t$ of metrics which connect a $U_0$-allowable metric to a
$U_1$-allowable metric, there is an open neighborhood $\NdZero$ of
$0\in\Targ$ and an integer $\theta_{\met_t}(\spinc)$ with the
property that for all generic $\eta_0\times\eta_1\in \NdZero$, we have
that
$$\theta_{\met_t,\eta_0\times\eta_1}(\spinc)=\theta_{\met_t}(\spinc).$$
Moreover, if $\met_t'$ is another path of metrics connecting a
$U_0$-allowable metric with a $U_1$-allowable metric, then
$$\theta_{\met_t}(\spinc)=\theta_{\met_t'}(\spinc).$$
\end{prop}

\begin{proof}
According to Sard's theorem, for any generic
$\eta_0\times\eta_1\in\Targ$, the fiber
$$M_{\met_t,\eta_0\times\eta_1}=
\{(s,t,D)\big| \Theta_{\met_s}(D)\in \Lag_0 + \eta_0, 
\Theta_{\met_t}(D)\in \Lag_1 + \eta_1, s\leq t\}$$
is a compact, canonically oriented,
zero-dimensional manifold. We investigate the conditions necessary to
show that the fiber misses the boundary of the domain of $\ParTheta$.

The set of points $\eta_0\times\eta_1\in\Targ$ for which the fiber
$M_{\met_t,\eta_0\times\eta_1}$ does not contain boundary points of
the form $(s,1,D)$ or $(0,t,D)$ is an open set which contains $0$
(since $\met_i$ is $U_i$-allowable for $i=0,1$). Let $\NdZero$ be a
connected neighborhood of $0$ in this set. Moreover, the fiber
$M_{\met_t,\eta_0\times\eta_1}$ cannot contain points of the form
$(t,t,D)$ if the spaces $\Lag_0+\eta_0$ and $\Lag_1+\eta_1$ are
disjoint; but $(\Lag_0+\eta_0)\cap (\Lag_1+\eta_1)\neq \emptyset$ is
equivalent to the condition that the image of $\delta(\eta_0-\eta_1)=0
\in H^2(Y;S^1)$. This is a codimension $b_1(Y)$ sub-torus
${\Wall}$ of $\Targ$, so its complement is dense.

Thus, for a dense set of perturbations $\eta_0\times\eta_1\in
\NdZero$, the fiber $M_{\eta_0\times\eta_1}$ is a smooth submanifold
which misses the boundary of the domain of $\ParTheta$. Moreover, given
two generic perturbations $\eta_0\times\eta_1, \eta_0'\times\eta_1'\in
\NdZero$, a generic path in $\NdZero$ misses the locus ${\Wall}$ as
well, since it has codimension $b_1(Y)>1$ (this is what distinguishes
the case where $b_1(Y)>1$ from the case $b_1(Y)=1$).  By Sard's
theorem, then, a generic such path induces a compact cobordism between
$M_{\eta_0\times\eta_1}$ and $M_{\eta_0'\times\eta_1'}$. By lifting to
the covering space $\OurSym$, one can easily see that the cobordism
respects the partitioning into $\SpinC$ structures. Thus, $\NdZero$
has the required property.

Fix $\met_t$ and $\met_t'$. Since the space of $U$-allowable metrics
is path connected and the space of metrics over $\Sigma$ is
simply-connected, we can connect $\met_t$ and $\met_t'$ by a
two-parameter family of metrics with $\metk_{t,0}=\met_t$,
$\metk_{t,1}=\met_t'$, $\metk_{0,t}$ is $U_0$-allowable and
$\metk_{1,t}$ is $U_1$-allowable. This, together with a small generic
perturbation, gives rise to a cobordism between
$M_{\eta_0\times\eta_1,\met_t}$ and $M_{\eta_0'\times\eta_1',\met_t'}$
(which once again respects the partitioning into $\SpinC$ structures).
This completes the proof of the proposition.
\end{proof}

The invariance statement in Proposition~\ref{prop:IndepOfMetBig} holds
when $b_1(Y)=1$ as well, but the argument is more involved. The point
is that ${\Wall}$ now separates $G$ into two components. The above
proof shows that if we pick two paths of metrics and two perturbations in
the same component, then the invariant remains unchanged. So, from now
on, we can drop the path $\met_t$ from the notation for $\theta$. We
must show, then, that the invariant is actually independent of the
component. It will be convenient to make use of the involution on the
set of $\SpinC$ structures introduced in
Section~\ref{sec:Introduction}, $\spinc\mapsto{\overline\spinc}$. This
is the map which sends the complex spinor bundle $W$ of $\spinc$ to
the same underlying real bundle, given its conjugate complex structure
(and naturally induced Clifford action). Note that this action fixes
those $\SpinC$ structures which arise from $\Spin$ structures, and
moreover $\Conj{\spinc+\ell}=\Conj{\spinc}-\ell$, for any
$\spinc\in\SpinC(Y)$ and $\ell\in H^2(Y;\Z)$. The invariant is
preserved by this involution in the following sense:

\begin{lemma}
\label{lemma:Involution}
Let $Y$ be an oriented three-manifold with $b_1(Y)>0$. Then,
$$\theta_{\eta_0\times \eta_1}(\spinc)=\theta_{-\eta_0\times
-\eta_1}({\overline \spinc})$$
\end{lemma}

\begin{proof}
To see this, we make use of the Serre duality map. Serre duality gives
rise to a natural involution on $\Jac$, with the property that
$\Conj{q+x}=\Conj{q}-x$ for any $q\in\Jac$ and $x\in
H^1(\Sigma;\R)$. For any given metric $\met_t$, this involution
preserves the theta divisor (by Serre duality), and fixes the points
associated to spin structures.  The involution on $\Jac$ can be lifted
to an involution of $\OurJac$ which we can assume fixes
$\Lag_0(\spinc)$. The involution then carries $\Lag_1(\spinc)$ to
$\Lag_1({\overline \spinc})$. Unfortunately, this involution is not
defined on the symmetric product.  Instead, it gives an involution on
the set of points in the symmetric product which map injectively to
the theta divisor, and more generally, it gives a relation. We write
$D_1\sim_{h}D_2$ if $\Theta_{h}(D_1)=\Conj{\Theta_{h}(D_2)}$.

For any $y\in[0,1]$, let $M_y$ denote the moduli space
$$M_y=\left\{(D,s,t)\Bigg| 
\begin{array}l
s\leq t, \\
\exists D', D'' \in \Sym^{g-1}(\Sigma),
D'\sim_{h_s} D, D'' \sim_{h_x} D, \\
\Theta_{h_s}(D')\in \Lag_0+\eta_0, \\
\Theta_{h_t}(D'')\in \Lag_1+\eta_1,
\end{array} \right\}, $$
where $x=s+y(t-s)$.  For generic $\eta_0$, $\eta_1$, $\met_t$, $M_y$
is a smooth cobordism from $M_0$ to $M_1$. This follows from the fact
that the set of points in $\Sym^{g-1}(\Sigma)$ which do not map
injectively into the theta divisor is empty if $g=2$ and has complex
codimension $2$ for generic $h$ if $g>2$
(see p. 250 of~\cite{ArbCornalGriffHar}).
As before, $M_y$ can be
partitioned according to $\SpinC$ structures; let $M_y(\spinc)$ denote
the corresponding set. The same dimesnion counting also shows that
for generic $\met_t$, $\eta_0$, and $\eta_1$,
\begin{eqnarray*}
M_0(\spinc)=M_{\met_t,\eta_0\times \eta_1}(\spinc)
&{\text{and}}&
M_1(\spinc)=M_{\met_t,-\eta_0\times -\eta_1}({\overline\spinc}).
\end{eqnarray*}
This proves the lemma.
\end{proof}

If $Y$ is an integral homology $S^1\times S^2$, the invariance of
$\theta$ is an easy consequence of this lemma:

\begin{prop}
\label{prop:IndepOfMetOneZ}
If $Y$ is an integral homology $S^1\times S^2$, then the invariant
$\theta(\spinc)$ is independent of metrics and perturbations as well.
\end{prop}

\begin{proof}
If $\spinc$ is not a torsion
class, then $\Lag_0(\spinc)$ and $\Lag_1(\spinc)$ are disjoint, so the
proof of Proposition~\ref{prop:IndepOfMetBig} applies.

For the torsion $\SpinC$ structure $\spinc$ (i.e. the one corresponding to the
spin structure), there are {\em a priori} two invariants, depending on
the sign $\delta(\eta_0-\eta_1)$ for the perturbation $\eta_0\times
\eta_1$. But Lemma~\ref{lemma:Involution} guarantees that
$$\theta_{\eta_0\times \eta_1}(\spinc)=\theta_{-\eta_0\times
-\eta_1}({\overline \spinc});$$ and since $\spinc$ comes from the spin
structure, ${\overline\spinc}=\spinc$, while
$$\delta(\eta_0-\eta_1)=-\delta(-\eta_0-(-\eta_1)),$$ so 
the invariants for both perturbations are equal.
\end{proof}

More generally, we have:

\begin{prop}
\label{prop:IndepOfMetOne}
If $Y$ is a rational homology $S^1\times S^2$, then the invariant
$\theta(\spinc)$ is independent of metrics and perturbations.
\end{prop}

This general case involves a study of spectral flow, and takes us
slightly out of the general framework of the preceding discussion (and
it is surprisingly more involved than the special case), so we
relegate that case to Section~\ref{sec:WallCross} (where it is
restated as Proposition~\ref{prop:IndepOfMetOneGen}).

Our final goal in this section is to prove
Proposition~\ref{prop:StandardProperties}.

\vskip0.3cm
\noindent{\bf{Proof of Proposition~\ref{prop:StandardProperties}}.}
The finiteness statement is clear from the fact that the fibers of
$\ParTheta$ are compact. 

The statement about involution invariance follows from
Lemma~\ref{lemma:Involution}. When we reverse the orientation of $Y$,
then the corresponding orientation of $H^*(Y)$ changes by a factor of
$(-1)^{b_1+1}$, so the last statement follows from the discussion
preceding Proposition~\ref{prop:IndepOfMetBig}. 
\qed
\vspace{0.2in}

In Section~\ref{sec:Splicing}, we establish the final step towards
proving topological invariance, showing that $\theta(\spinc)$ remains
invariant under stabilization. It must be shown that if one begins
with $\Sigma$, then the associated invariant agrees with the invariant
calculated from the Riemann surface obtained by a connected sum
$\Sigma\# (S^1\times S^1)$, using the handlebodies $U_0\# (D^2\times
S^1)$ and $U_1 \# (S^1\times D^2)$. This result fits naturally into
the context of a splicing construction, which also proves a technical
device which underpins the calculations of the invariants. We defer
these results to Section~\ref{sec:Splicing}, addressing first the
lemmas which we have stated thus far without proof.

\section{Degenerating metrics}
\label{sec:MoveTheta}

In the definition of $\theta(\spinc)$, we used several facts
(Lemmas~\ref{lemma:MissTheta} and \ref{lemma:MissThetaHS}) about the
behaviour of the theta divisor under certain degenerations of the
metric $\met$ over $\Sigma$. The aim of this section is to provide
proofs of these results. We opt to give analytical proofs, which rely
only on elementary properties of the $\DBar$ operator on cylindrical
manifolds. The statements given here, though, can be interpreted in
terms of algebraic geometry, where they address degenerations of the
theta divisor in a family of curves acquiring nodal singularities. For
a related discussion from this point view, see~\cite{Fay}.  (This
remark pertains also to the discussion in Section~\ref{sec:Splicing}.)

The proof of both Lemmas~\ref{lemma:MissTheta} and
\ref{lemma:MissThetaHS} involve local compactness and then passing to
cylindrical end models, as is familiar in gauge theory (though the
results here are considerably more elementary than is typical in gauge
theory).  More precisely, thanks to local compactness, points in the
theta divisors for Riemann surfaces undergoing suitable degenerations
give rise to points in the $L^2$ theta divisor for the cylindrical
end model of $\Sigma-(\gamma_1\cup...\cup\gamma_n)$, where the
$\gamma_i$ are embedded, disjoint, closed curves in $\Sigma$.  We then
appeal to basic results about the $L^2$ kernel of the Dirac operator.
To state these results, we set up some notation.

Let $\FComp$ be a compact, oriented two-manifold with $\NumEnds$ boundary
circles, given a product metric in a neighborhood of its boundary. Let
$$
\FCyl= \FComp\cup_{\partial \FComp}\coprod_{j=1}^{\NumEnds} (S^1_{(j)}\times
[0,\infty)) $$ be the associated complete manifold
with cylindrical ends, and $$ \FFill = \FComp\cup_{\partial
\FComp}\coprod_{j=1}^{\NumEnds}
\CDisk_{(j)} $$ be the associated compact Riemann surface; here,
$\CDisk_{(j)}$ is a copy of the two-dimensional disk with a product
metric near its boundary. Conformally, $\FFill$ is obtained from
$\FCyl$ by adding $\NumEnds$ points ``at infinity'' $\{p_1,...,p_\NumEnds\}$,
which correspond to the centers of the attached disks in the
description of $\FFill$. Fix a spin structure ${\spin}_0$ on $\FFill$
with associated spinor bundle $E$ throughout the following discussion
(for instance, when $\FComp$ has genus zero, this spin structure is
uniquely determined). This canonically induces a spin structure on
$\FCyl$, which gives rise to a canonical connection $B_0$ on the
spinor bundle $\ECyl$ over $\FCyl$ with normalized curvature form.  Up
to gauge transformations, any other connection with normalized
curvature form differs from $B_0$ by a cohomology class in
$H^1(\FCyl;\R)$.

Now, let $B$ be any connection on $\ECyl$ over $\FCyl$ with normalized
curvature form.  We can relate the $L^2$-extended $B$-harmonic spinors
on $\FCyl$ with holomorphic data on $E$ over
$F$. Following~\cite{APSI}, a $B$-harmonic spinor $\Psi$ is said to be
an $L^2$-extended section if it is in $L^2_{\loc}$, and over each
cylinder $S^1_{(j)}\times [0,\infty)$, we can write the restriction of
$\Psi$ as a sum of an $L^2$ $B$-harmonic section plus a constant
(hence $B$-harmonic) section. The holomorphic data over $F$ is
obtained by viewing $E$ as a holomorphic vector bundle over $F$ (by
using the $\DBar$-operator associated to the spin-connection coming
from $\spin_0$).  Given $p\in F$, let $\Ideal{p}$ denote the ideal
sheaf at $p$.  We have the following result (a related discussion can
be found in~\cite{APSI}):

\begin{prop}
\label{prop:LTwoCohom}
Let $B$ be a connection on $\ECyl$ with normalized curvature form, and write
$B=B_0+i\xi$ for $\xi\in H^1(\ECyl;\R)$. Let
$$\xi_j=\langle \xi, [S^1_{(j)}]\rangle.$$
Then, the space of $L^2$-extended
holomorphic sections of $\ECyl$ over $\FCyl$ is canonically identified with
$$H^0(\FFill,\EFill\otimes
\Ideal{p_1}^{\lceil \xi_1 -  \OneHalf \rceil}\otimes
... \otimes\Ideal{p_{\NumEnds}}^{\lceil \xi_{\NumEnds} - \OneHalf
\rceil}
).$$
Here, $\lceil x \rceil$ denotes the smallest integer greater than $x$.
\end{prop}

\begin{proof}
Note that the cylindrical-end metric $\MetCyl$ on $\FCyl$ is conformal
equivalent to the metric $\MetComp$ on $\FFill-\{p_1,...,p_\NumEnds\}$
inherited from $\FFill$. Indeed, we can write $$\MetCyl =
e^{2\tau}\MetComp,$$ where $$\tau\colon \FCyl\longrightarrow
[0,\infty)$$ is a real-valued function, which agrees with the real
coordinate projection on the cylindrical ends. The spinor bundles of
the two manifolds can be (metrically) identified accordingly, with a
change in the Clifford action to reflect the conformal change. With
respect to this change, (see~\cite{Hitchin}, bearing in mind that we
are in two dimensions), the Dirac operator over $\FCyl$ can be written:
$$\DiracCyl = e^{-\OneHalf \tau}\DiracDisk e^{\OneHalf \tau};$$ i.e.,
multiplication by $e^{-\OneHalf\tau}$ induces a vector space
isomorphism from the $\DiracDisk$-harmonic spinors to the
$\DiracCyl$-harmonic ones. Moreover, a section of a bundle is in $L^2$
for the cylinder iff $e^{\tau}\phi$ is in $L^2$ for
$\FFill-\{p_1,...,p_\NumEnds\}$.

From this discussion, it follows that the $L^2$-harmonic spinors on
$\FCyl$ are identified with the space of harmonic spinors over
$\FFill-\{p_1,...,p_\NumEnds\}$ for which $e^{\frac{\tau}{2}}\phi$
lies in $L^2$ (for $\FFill-\{p_1,...,p_\NumEnds\}$). The 
proposition follows from this, along with some considerations in the
neighborhoods of the punctures.

Consider $\CDisk-\{0\}$, with the trivial line bundle endowed with a
connection $B=d+i\xi d\theta$ -- this is the model of the punctured
neighborhood of the $p_j\in \FFill$.  Under the standard
identification $S^1\times [0,\infty)\cong \CDisk-\{0\}$, the function
$e^{-t}$ (which is $e^{-\tau}$ over $\FCyl$) corresponds to the
radial coordinate $r$ on the disk. Moreover, multiplication by $e^{\xi
t}=r^{-\xi}$ induces an isomorphism from the space of (ordinary)
holomorphic functions on $\CDisk-\{0\}$ to the space of
$\DBar_B$-holomorphic functions. Under these correspondences, a
holomorphic function $\phi$ which vanishes to order $k$ corresponds to a
$L^2$-harmonic spinor on $\FCyl$ iff
\begin{eqnarray*}
\int |\phi r^{-\xi- \OneHalf}|^2 r dr d\theta  &\leq &
C \int r^{2k-2\xi}dr d\theta
\end{eqnarray*}
is bounded, i.e. iff $k> \xi-\OneHalf$. In the borderline case where
$\xi\equiv \OneHalf$, the holomorphic functions on $\CDisk-\{0\}$
which vanish to order $\xi-\OneHalf$ correspond to harmonic sections
over the cylinder whose pointwise norm is bounded, and hence they lie
in the space $L^2$ (of the cylinder) extended by constants.

The proposition follows.
\end{proof} 

Given this proposition, then, we can give a proof of
Lemma~\ref{lemma:MissTheta}.

\vskip0.3cm
\noindent
{\bf{Proof of Lemma~\ref{lemma:MissTheta}.}}  Suppose that
$\{\met_i\}_{i=1}^{\infty}$ is a sequence of metrics whose
neck-lengths along the $\gamma_i$ all go to infinity, and whose
restrictions away from the necks lie in a compact family of metrics on
the genus zero surface $F^c=\Sigma-\nbd{\gamma_1\cup...\cup\gamma_g}$.
Let $B_i$ be a sequence of connections which lie in the theta divisor
of $\met_i$. This means that we can find a sequence of non-zero
sections $\phi_i$ of $E$ over $\Sigma$ (with metric $\met_i$), so that
$\DBar_{B_i}\phi_i=0$. By renormalizing, we can assume without loss of
generality that the $\sup_{\Sigma}|\phi_i|=1$. Since the metric in a
neighborhood of the tubes is flat, the supremum is always achieved in
the compact piece $\FComp\subset \Sigma$.  After passing to a
subsequence, the $B_i$ converge (locally in $\Cinfty$) to a connection
$B_\infty$ on $\FCyl$ with normalized curvature form. In fact, local
compactness of holomorphic functions ensures that (after passing to a
subsequence) the $\phi_i$ converge locally (in $\Cinfty$) to a
$B_{\infty}$-holomorphic section $\phi_\infty$. Once again, the
supremum of $|\phi_\infty|$ must be $1$, so in particular,
$\phi_\infty$ is a non-vanishing, $L^2$-extended,
$B_{\infty}$-holomorphic section.  From
Proposition~\ref{prop:LTwoCohom}, it then follows that the holonomy of
$B_{\infty}$ around at least one of the boundary circles must be
congruent to $\OneHalf$ (modulo $\Z$). Specifically, the spin
structure $E_0$ on $F=S^2$ has degree $-1$, so the dimension of the space
of $\Sobol{2}{}$-bounded, holomorphic sections is calculated by the
formula $$\dim_{\C}H^0(\FFill,\EFill_0\otimes
\Ideal{p_1}^{\lceil \xi_1 -  \OneHalf \rceil}\otimes
... \otimes\Ideal{p_{2g}}^{\lceil \xi_{2g} - \OneHalf
\rceil})= \min\left(0, \sum_{i=1}^{2g} \lceil \xi_i - \OneHalf \rceil\right).$$
Note that the holonomies $\xi_i$ all add up to
zero, and, since the
ends of $\FComp$ are naturally come in pairs $S^1_{(i)}$, 
$S^1_{(i+g)}$ with $\xi_i \equiv -\xi_{i+g} \pmod\Z$, it
follows that the dimension is non-zero only if at least one of the
holonomies is $\OneHalf$ modulo $\Z$.

Strictly speaking, to apply Proposition~\ref{prop:LTwoCohom}, we note
that natural compactification of $F^c$ is a sphere, and the two
methods for measuring holonomy -- comparing holonomies against any
spin structure which extends over $U$ (which is used in the statement
of Lemma~\ref{lemma:MissTheta}) and comparing against the spin
structure which extends over the sphere (which we use in the statement
of Proposition~\ref{prop:LTwoCohom}) -- coincide. This is obvious from
the Kirby calculus description of $U$.
\qed
\vspace{0.2in}

The proof of Lemma~\ref{lemma:MissThetaHS} is analogous.

\vskip0.3cm
\noindent
{\bf{Proof of Lemma~\ref{lemma:MissThetaHS}.}} As the metric is
stretched in a sequence $\met_i$, any sequence of points $B_i\in
\Theta_{h}(\Sym^{g-1}(\Sigma))\cap \Hol_{\gamma_1}^{-1}(0)\cap...\cap
\Hol_{\gamma_{g-1}}^{-1}(0)$ has a subsequence which converges to a
connection $B_{\infty}$, which now can be viewed as a connection on a
torus with cylindrical ends $\FCyl$.  Kernel elements then converge to
a section which, according to Proposition~\ref{prop:LTwoCohom}
corresponds naturally to a holomorphic section of a line bundle
$\EFill$ over a compact torus $\FFill$. But there is only one spin
structure on $\FFill$ which admits harmonic spinors (the trivial
bundle), and it corresponds to the spin structure on $\FFill$ which
does {\em not} bound. Thus, around any curve which bounds in $U$, the
difference in the holonomy between this spin structure and any spin
structure which bounds in $U$ is $1/2$.
\qed \vspace{0.2in}

\section{Splicing}
\label{sec:Splicing}

The aim of this section is to give a more detailed analysis of the
theta divisor, using a splicing construction whose consequences
include the stabilization invariance of $\theta(\spinc)$, and a
technical result which will be of importance in subsequent sections.

We introduce notation.  Fix a pair of Riemann surfaces $\Sigma_i$ for
$i=1,2$, and let $\Sigma_i^c$ be the complement in $\Sigma_i$ of an
open disk centered at $p_i\in \Sigma_i$, endowed with a product-end
metric (which we can extend over the disks to obtain the metrics over
$\Sigma_i$).  Let $\Sigma_i^+$ (for $i=1,2$) denote cylindrical-end
models of the surfaces, $$\Sigma_i^+=\Sigma_i^c\cup
\big(S^1\times [0,\infty)\big).$$ Let $\Sigma_1\#_T\Sigma_2$ denote
the model for the connected sum 
of $\Sigma_1$ and $\Sigma_2$, a surface of genus $g=g_1+g_2$
with a neck length of $2T$; i.e.
$$
\Sigma_1\#_T\Sigma_2=\Sigma_1^c\cup \big([-T,T]\times S^1\big) \cup
\Sigma_2^c.$$ 

Fix non-negative integers $k_1$, $k_2$ so that $k_1+k_2 = g_1 + g_2 -
1$. For all $T$, there is an obvious natural map $$\Splice{T}\colon
\Sym^{k_1}(\Sigma_1^c)\times \Sym^{k_2}(\Sigma_2^c) \longrightarrow
\Sym^{g-1}(\Sigma_1\#_T\Sigma_2).$$ Fix a spin structure $\spin_0$
over $\Sigma_1\#\Sigma_2$. This allows us to compare the various
Jacobians as $T$ varies; i.e. it gives us identifications:
$$\Jac(\Sigma_1\#_T\Sigma_2)\cong H^1(\Sigma_1\#_T\Sigma_2;S^1).$$
Similarly, we can use the natural extension of this structure over
$\Sigma_1$, $\Sigma_2$ to fix identifications
\begin{eqnarray*}
\Jac(\Sigma_1)\cong H^1(\Sigma_1;S^1)&{\text{and}}&
\Jac(\Sigma_2)\cong H^1(\Sigma_2;S^1).
\end{eqnarray*}

To state the splicing result, we use Abel-Jacobi map, thought of as
follows. Fix a Riemann surface $\Sigma$ with a basepoint $p$, then
$$\AJ^{(k)}\colon\Sym^k(\Sigma)\longrightarrow\Jac(\Sigma)$$ is the map
which takes an effective divisor $D\in\Sym^k(\Sigma)$ to the unique
connection $A$ with normalized curvature form, which admits a
$\DBar_A$-meromorphic section $\phi$ whose associated divisor is
$D+(g-1-k)p$. When $k=g-1$, then we do not need a base point, and
$\AJ^{(g-1)}$ agrees with the map $\Theta$ from
Section~\ref{sec:Introduction}.

\begin{theorem}
\label{thm:Splicing}
In regions of the symmetric products supported away from the points
$p_i$, the composite of $\Splice{T}$ with
$\Theta_{\Sigma_1\#\Sigma_2}$ is homotopic to the product of
Abel-Jacobi maps.  Indeed, for any non-negative $k_1, k_2$ with
$k_1+k_2=g-1$, we have that the composite $$\begin{CD}
\Sym^{k_1}(\Sigma_1^c)\times \Sym^{k_2}(\Sigma_2^c)
@>{\Splice{T}}>>
\Sym^{g-1}(\Sigma_1\#_T\Sigma_2) 
@>{\Theta_{\Sigma_1\#\Sigma_2}}>>
H^1(\Sigma_1\#\Sigma_2;S^1)
\end{CD}$$
converges in the $C^1$ topology, as $T\goesto\infty$, to the map 
$$\begin{CD}
\Sym^{k_1}(\Sigma_1^c)\times \Sym^{k_2}(\Sigma_2^c)
@>\AJ_1\times\AJ_2>> H^1(\Sigma_1;S^1)\times
H^1(\Sigma_2;S^1)\cong H^1(\Sigma_1\#\Sigma_2;S^1),
\end{CD}$$
where $\AJ_i$ denotes the Abel-Jacobi map with basepoint $p_i$
$$\AJ_i=\AJ_i^{(k_i)}\colon \Sym^{k_i}(\Sigma_i)\longrightarrow
H^1(\Sigma_i;S^1).$$
\end{theorem}

\begin{remark}
The seasoned gauge theorist will identify the last vestiges of a
``gluing theorem'' here. However, the present result is significantly
easier than the usual gluing results.
\end{remark}

Before giving the proof of Theorem~\ref{thm:Splicing}, we recall the
construction of the map $$\Theta_{\met}\colon
\Sym^{g-1}(\Sigma)\rightarrow H^1(\Sigma;S^1).$$
For a given divisor $D\in\Sym^{g-1}(\Sigma)$, $\Theta_{\met}(D)$ is
the unique connection $B$ in the spinor bundle $E$ with normalized
curvature form which admits a $\DBar_B$-holomorphic section whose
associated divisor is $D$. To find it, first fix a section $\phi$ of
$E$ whose vanishing set is $D$; then find any connection $A$ on $E$
for which $\phi$ is $\DBar_A$-holomorphic.  Now, let $f$ be a function
which solves $$ id*df=F_{B_0}-F_A.  $$ In the above equation (and
indeed throughout this section), $B_0$ denotes the connection with
normalized curvature form on the spinor bundle induced by the spin
structure $\spinc_0$. Then, $A+ i *df$ will represent
$\Theta_{\met}(D)$.  The key to Theorem~\ref{thm:Splicing}, then, is
to select the initial connection $A$ carefully. Before giving the
proof, we name one of the fundamental objects which arises in the
construction.

\begin{defn}
If $D$ is a divisor of degree $g-1$ and $(B,\phi)$ is a connection
with normalized curvature form for which $\DBar_B\phi=0$, and
$\phi^{-1}(0)=D$, then we call $(B,\phi)$ a holomorphic pair
representative for the divisor $D$. Of course, the gauge equivalence
class of $B$ represents $\Theta_{\met}(D)$.  
\end{defn}

\vskip0.3cm
\noindent{\bf{Proof of Theorem~\ref{thm:Splicing}.}}
Since one of the $k_i>g_i-1$, we can assume without loss of generality
that $k_2>g_2-1$. 
Pick a partition of unity $\psi_1,\psi_2$ over $[-2,2]\times S^1$
subordinate to the cover
$$\big\{[-2,1)\times S^1, (-1,2]\times S^1 \big\}.$$
We can transfer this partition of unity (by extending by constants in
the obvious way) to a partition of unity on
$\Sigma_1\#_T\Sigma_2$
subordinate to the cover
$$\big\{\Sigma_1^c\cup([-T,1)\times S^1), ((-1,T]\times S^1)\cup
\Sigma_2^c\big\}$$
(provided that $T>2$). We denote this partition of unity also by
$\{\psi_1,\psi_2\}$ although, technically, it does depend on
$T$. However, notice that for all $T$, the ($L^2$ and $\Cinfty$) norms
of $d\psi_i$ remain constant.  

Fix a pair of divisors $D_i\in
\Sym^{k_i}(\Sigma_i)$.
After deleting the points $p_1$ and $p_2$, we find suitably normalized
holomorphic pair representatives for the $D_i$ over the
cylindrical-end manifolds; i.e.  connections $A_1$, $A_2$ with fixed
curvature form on $\Sigma_1^+$, $\Sigma_2^+$ and sections $\phi_1,
\phi_2$ whose vanishing locus is $D_1, D_2$ respectively, with
asymptotic expansions (with respect to some trivialization of the
spinor bundle over the flat cylinder) of the form:
\begin{eqnarray}
\phi_1= e^{\alpha t} + \BigO(e^{(\alpha-1)t})&{\text{and}}&
\phi_2= e^{-\alpha t} + \BigO(e^{(-\alpha-1)t}),
\label{eq:Decay}
\end{eqnarray}
where $\alpha={-g_1+k_1+\OneHalf}= g_2-k_2-\OneHalf$. Note that the
leading terms in the asympotic expansions here have a particularly
simple form; this can be arranged by first untwisting the imaginary
part of the leading term using a gauge transformation, and then by
rescaling the $\phi_i$ by real constants if necessary. Note that the
decay rates come from Proposition~\ref{prop:LTwoCohom}: according to
that proposition, sections with the prescribed decay for $\phi_i$
correspond to sections of the spinor bundles over $\Sigma_i$ which
vanish to order $-g_i+k_i+1$ at the connected sum point.  Starting
from these sections $\phi_1$, $\phi_2$, we will construct for all $T$
holomorphic pair representatives $A_1\#_{T}A_2$ for
$\Splice{T}(D_1,D_2)$, and show that the gauge equivalence classes of
$A_1\#_{T}A_2$ converge, as $T\goesto\infty$ to those of $A_1$ and
$A_2$.

There is a natural connection with normalized curvature form
on $\Sigma_1\#_T\Sigma_2$ induced from
$A_1$ and $A_2$, which we write as $A_1\#A_2$, which is obtained from
the identification
$\Jac(\Sigma_1\times\Sigma_2)\cong H^1(\Sigma_1\#\Sigma_2; S^1)\cong
H^1(\Sigma_1; S^1)\times H^1(\Sigma_2; S^1)\cong
\Jac(\Sigma_1)\times\Jac(\Sigma_2)$ coming from our fixed spin
structure $\spin_0$. Consider, then, the section $\phi_T=\psi_1
\phi_1 + e^{2\alpha T}\psi_2\phi_2$. For all sufficiently large $T$,
this section does not vanish in the neck region $[-1,1]\times S^1$
(indeed, the restriction of $\phi_T$ to this region is $e^{\alpha
(T+t)}+\BigO(e^{(\alpha-1)T})$). Although $\phi_T$ is not
$\DBar_{A_1\#A_2}$ holomorphic, it is holomorphic for the
$\DBar$-operator 
$$
	\DBar_{A_1\# A_2} 
	   - (\DBar\psi_1)\frac{\phi_1}{\phi_T}
	   - e^{2\alpha T} (\DBar \psi_2)\frac{\phi_2}{\phi_T};
$$ 
so if we let $\Error$ be the form 
\begin{equation}
\label{eq:Error}
	\Error=
		\mathrm{Im}\left((\DBar\psi_1)\frac{\phi_1}{\phi_T} +
		e^{2\alpha T} (\DBar \psi_2) \frac{\phi_2}{\phi_T}\right),
\end{equation}
then $\phi_T$ is
holomorphic for $$A_3= A_1\# A_2 + 2i\Error.$$ The connection $A_3$ is a
good first approximation to the desired connection corresponding to
$\Splice{T}$.

The curvature form is normalized once we find $f$ so that
$d*df=-2d\Error$. We would like to show that this does not change the
cohomology class by much; i.e. as the tube length $T$ is increased,
the cohomology correction $2\Error+*df$ tends to zero (viewed as
elements in $H^1(\Sigma^c_1\coprod \Sigma^c_2;\R)\cong
H^1(\Sigma_1\#\Sigma_2;\R)\cong H^1(\Sigma^+_1;\R)\oplus
H^1(\Sigma^+_2;\R)$). To do this, we find it convenient to use harmonic
forms.  Let $\Harm_T$ denote the space of harmonic one-forms on
$\Sigma_1\#_T\Sigma_2$, $$\Harm_T=\{a\in
\Omega^1(\Sigma_1\#_T\Sigma_2)\big| da=d*a=0\},$$ and let
$\Proj_{T}$ denote the $\Sobol{2}{}$-projection to $\Harm_T$. By Hodge
theory, the map from closed one-forms, given by $z\mapsto
[\Proj_T(z)]$, induces the identity map in cohomology. Moreover, it is
identically zero on co-closed one-forms (by
integration-by-parts). Thus,
$$[\Proj_T(2\Error+*df)]=[\Proj_T(2\Error)].$$ Note that
$\lim_{T\goesto\infty}\|\Error\|_{\Sobol{2}{}}=0$ (indeed, using the
fact that $\DBar(\psi_1+\psi_2)=0$, the decay condition in
Equation~\eqref{eq:Decay} and the expression for $\Error$,
Equation~\eqref{eq:Error}, it is easy to see that
$\|\Error\|_{\Sobol{2}{}}=\BigO(e^{(\alpha-1)T})$), so
$\lim_{T\goesto\infty}\Proj_T(\Error)=0$ in $L^2$. The fact that the
harmonic projections tend to zero, then, is a consequence of elliptic
regularity, as follows:

\begin{lemma}
Let $h_T$ be a sequence of harmonic forms on $\Sigma_1\#_T\Sigma_2$ whose
$\Sobol{2}{}$ norm tends to zero. Then the cohomology classes $[h_T]$
tend to zero, as well.
\end{lemma}

\begin{proof}
Since the operator $d+d^*$ is translationally invariant in the
cylinder $S^1\times\R$, there is a single constant $C$ which works for
all the manifolds $\Sigma_1\#_T\Sigma_2$, so that for any form $\phi\in\Wedge^*(\Sigma_1\#_T\Sigma_2)$,
$$\|\phi\|_{\Sobol{2}{k+1}}
\leq C(\|\phi\|_{\Sobol{2}{}} + 
\|(d+d^*)\phi\|_{\Sobol{2}{k}})$$
(see for instance~\cite{FreedUhlenbeck}); thus,
$$\|h_T\|_{\Sobol{2}{k}}\leq C(\|h_T\|_{\Sobol{2}{}}).$$ This together
with the Sobolev lemmas shows that the forms $h_T$ converge to zero in
$\Cinfty$ over any compact set. But the cohomology class of any of the
$h_T$ is determined by its restriction to the (compact) subset
$\Sigma_1^c\coprod\Sigma_2^c\subset
\Sigma_1\#_T\Sigma_2$. 
\end{proof}

This proves the convergence of $\Splice{T}$ in $C^0$. To prove $C^1$
convergence, we argue that any path in the space of divisors over
$\Sigma_1^c$ and $\Sigma_2^c$ can be covered by a path in the space of
holomorphic pairs whose derivatives satisfy decay conditions analogous
to Equation~\eqref{eq:Decay}. To see this, it helps to consider
Fredholm deformation theory for the symmetric product which arises by
viewing the latter space as the zeros of a non-linear equation on the
cylinder. More precisely, let $Z_i$ be a vector space of
compactly-supported forms in $\Sigma^c_i$ which map isomorphically to
$H^1(\Sigma^c_i)$. Consider the the map
$$\Omega^{0,0}_{\delta_i}(E_i^+)\times Z_i\longrightarrow
\Omega^{0,1}_{\delta_i}(E_i^+),$$ given by
$\DBar_{B_0}\Phi+(ia)^{0,1}\Phi$, with weighted Sobolev topologies on
the $\Omega^{0,*}(E_i^+)$, i.e.
$$\|\Phi\|_{\delta_i}=\|e^{{-\frac{\delta_i
t}{2}}}\Phi\|_{\Sobol{2}{}},$$ where $\delta_i=g_i-1-k_i$, and $t$ is
a smooth function on $\Sigma_i^+$ which extends the real coordinate
function on the cylinder $S^1\times[0,\infty)$ (the decay rate is
chosen to for $\Phi$ to correspond to a holomorphic section of the
degree $k_i$ bundle over $\Sigma_i$, according to
Proposition~\ref{prop:LTwoCohom}). This is a non-linear, Fredholm map
whose zero locus (away from the trivial $\Phi\equiv 0$ solutions) is
transversally cut out by the equations. This zero locus, the space of
holomorphic pairs on $\Sigma_i^+$, admits a natural submersion to the
symmetric product, given by taking $(A,\Phi)$ to the divisor where
$\Phi$ vanishes (this models the quotient by the natural $\C^*$ action
on the space of holomorphic pairs). Thus, any tangent vector in
$\Sym^{k_i}(\Sigma_i^c)$ can be represented by a pair $(a,\phi)\in
Z_i\times \Omega^{0,0}(E_i^+)$ where $a$ is compactly supported in
$\Sigma_i^c$ and $$\phi=C
e^{(-g_i+k_i+\OneHalf)t}+\BigO(e^{(-g_i+k_i-\OneHalf)t}),$$ for some
constant $C$ (depending on the tangent vector).

Now, consider a pair of smooth paths $D_1(s)$ and $D_2(s)$, and a
corresponding paths of holomorphic pairs $(A_1(s),\phi_1(s))$ and
$(A_1(s),\phi_2(s))$.  Note that the derivative $\frac{d}{ds}
(A_1(s)\#_T A_2(s))$, restricted to $\Sigma_1^c\coprod \Sigma_2^c$, is
the differential of $\Theta_{\Sigma_1}\times \Theta_{\Sigma_2}$. To
prove the $C^1$ convergence, we must show that the derviative of the
error term converges to zero, i.e. writing 
$$
	[\frac{d}{ds} A_3(s)] =
		[\frac{d}{ds} A_1(s)\# A_2(s)] 
		+ 2i[\Proj_T(\frac{d}{ds} \Error(s))],
$$
we must show that $[\Proj_T(\frac{d}{ds} \Error)]\goesto 0$ in $T$.
Note first that
\begin{eqnarray*}
\frac{d}{ds}\Error(s)&=& (\DBar \psi_1)\left(\frac{d \phi_1}{ds}\frac{1}{\phi_T} -
\frac{d \phi_T}{ds}\frac{\phi_1}{\phi_T^2}\right) + 
e^{2\alpha T}(\DBar \psi_2)\left(\frac{d \phi_2}{ds}\frac{1}{\phi_T^2} -
\frac{d \phi_T}{ds}\frac{\phi_2}{\phi_T^2}\right) \\ 
&=& (\DBar \psi_1)\left(\frac{d \phi_1}{ds}\frac{1}{\phi_T} -
\left(\psi_1 \frac{d \phi_1}{ds}+
e^{2\alpha T} \psi_2 \frac{d \phi_2}{ds}\right) 
\frac{\phi_1}{\phi_T^2}\right) \\
&& + 
e^{2\alpha T}(\DBar \psi_2)\left(\frac{d \phi_2}{ds}\frac{1}{\phi_T} -
\frac{d \phi_T}{ds}\frac{\phi_2}{\phi_T^2}\right).
\end{eqnarray*}
It is easy to see from this 
that the derivative of the error is supported in the region
$[-1,1]\times S^1\subset \Sigma_1\#_T\Sigma_2$, and it is universally
bounded (e.g. in $C^0$) independent of $T$. For example, since
\begin{eqnarray*}
\phi=
e^{\alpha (t+T)}+ \BigO(e^{(\alpha-1) T}) &{\text{and}}&
\frac{d\phi_2}{ds}
= C e^{-\alpha (t+T)} + \BigO(e^{(-\alpha-1)T})
\end{eqnarray*} 
for $t\in [-1,1]\times S^1$ (according to
Fredholm perturbation theory we discussed above), we see that
$|e^{\alpha T}(\DBar \psi_2)\frac{d\phi_2}{ds}\frac{1}{\phi_T}|$ is
bounded above for all $T$. (The other terms follow in a similar
manner.)  Since, moreover, there is a universal constant $K$
independent of $T$ so that for any harmonic form $h$, $$\|h(0)\|\leq K
e^{-T} (\|h(-T)\| + \|h(T)\|)$$ (this follows from standard asymptotic
expansion arguments see for instance~\cite{APSI},
\cite{Yoshida}), the harmoic projection
$\Proj_T(\frac{d}{ds}\epsilon(s))$ converges to zero exponentially.
This proves the $C^1$ convergence statement, and concludes the proof
of Theorem~\ref{thm:Splicing}.
\qed
\vspace{0.2in}

Theorem~\ref{thm:Splicing} can be used to prove the
``stabilization invariance'' of the
invariant we are studying.  Specifically, in
Section~\ref{sec:DefTheta}, we gave a definition of the invariant
$\theta(\spinc)$, which refers to a choice of Heegaard decomposition
for $Y$. In the next proposition, we show that it is independent of that,
depending only on the underlying three-manifold.

\begin{prop}
\label{prop:Stabilization}
The invariant $\theta(\spinc)$ is independent of the Heegaard
decomposition used in its definition, thus it gives a well-defined
topological invariant.
\end{prop}

\begin{proof}
Fix a genus $g$ Heegaard decomposition of $Y=U_0\cup_{\Sigma}
U_1$. There is a ``stabilized'' genus $g+1$ Heegaard decomposition of
$Y$, corresponding to the natural decomposition $$S^3=(S^1\times\CDisk)\cup_{S^1\times S^1} (\CDisk\times S^1);$$ i.e. let
\begin{eqnarray*}
U_0'= U_0 \# (S^1 \times \CDisk), \Sigma'=\Sigma\#(S^1\times S^1) &{\text{and}}&
U_1'= U_1 \# (\CDisk\times S^1),
\end{eqnarray*}
and consider the Heegaard decomposition $$Y=U_0'\cup_{\Sigma'}U_1'.$$
We would like to show that the invariant $\theta$ associated to the
Heegaard decomposition $U_0\cup_{\Sigma}U_1$ agrees with that
associated to $U_0'\cup_{\Sigma'}U_1$, which we will denote
$\theta'$. 

Fix a metric on the torus $S^1\times S^1$.  We observe that one can
find $U$-allowable metrics $\met$ on a Riemann surface $\Sigma$ with
the property that for all sufficiently large $T$, $\met\#_T(S^1\times
S^1)$ are $U'$-allowable. To see this, let $\met(T_1)$ denote the
metric on $\Sigma$ which is stretched out along $g$ of the attaching
circles of $\Sigma$. We show there is a $T_0$ so that for all
$T_1,T_2>T_0$, $\met(T_1)\#_{T_2}(S^1\times S^1)$ is
$U'$-allowable. If this were not the case, we could stretch both
tube-lengths simultaneously, and extract a subsequence of spinors, which
would converge to a non-zero harmonic spinor either on the punctured
$(S^1\times S^1)$ (with a cylindrical end attached) -- which cannot
exist in light of the holonomy constraint coming from $S^1\times
\CDisk$, see Proposition~\ref{prop:LTwoCohom} -- or a harmonic
spinor on the genus zero surface with $g+1$ punctures obtained by
degenerating the punctured version of $\Sigma$. This is ruled out by
the holonomy constraints at infinity, as in the proof of
Lemma~\ref{lemma:MissThetaHS} (the holonomy around the curves
corresponding to the attaching circles vanish as in the proof of that
lemma; around the curve corresponding to the connected sum neck it
vanishes since that curve bounds in $\Sigma'$). Thus, for sufficiently
large $T_1$, the metric $\met(T_1)$ has the desired properties.

In view of this observation, we can find a path of metrics $\met_t$ on
$\Sigma$ to calculate $\theta$, with the property that for all
sufficiently long connected sum tubes, the family of metrics obtained
by connecting $\met_t$ with a constant metric on the torus
$F=S^1\times S^1$ can be used to calculate $\theta'$.  Choose a point
$p\in\Sigma$. The fiber of the map $$\ParTheta\colon
\Sym^{g-1}(\Sigma)\times [0,1]\times [0,1]\longrightarrow \Targ$$
(used in definition of $\theta$) over a generic point
$\eta\in\Targ$ misses the submanifold of divisors
$D\in\Sym^{g-1}(\Sigma)$ which contain the point $p$. In other words,
there is a compact region $K\subset
\Sigma-p$ so that $\ParTheta^{-1}(\eta)\subset
\Sym^{g-1}(K)\times[0,1]\times[0,1])$. 

Consider the one-parameter family of maps
$$\ParTheta'_T\colon\Sym^{g-1}(\Sigma\#_T F)
\longrightarrow\TargP,$$
used in defining the invariant $\theta'$ for the Heegaard
decomposition $U_0'\cup_{\Sigma'}U_1'$ (using metrics with length
parameter $T$).  Note that we have an isomorphism $$\TargP\cong
\Targ\times \Jac(F).$$ Under this isomorphism, the origin corresponds
to $0\times \spinc_0$, where $\spinc_0$ is a spin structure on the
torus which bounds. Let $q\in F$ be the pre-image of $\spinc_0$ under
the Abel-Jacobi map 
$$\AJ^{(1)}\colon F\longrightarrow \Jac(F).$$ Since
$\spinc_0$ admits no harmonic spinors, it follows that $q\in F$ is not
the connected sum point. Given a sequence of points $(D_T,s_T,t_T)\in
{\ParTheta'_T}^{-1}(\eta\times 0)$ with $T\goesto\infty$ using
compactness on the $\Sigma$-side, we obtain a subsequence which
converges to a divisor $D\in
\Sym^{g-1}(\Sigma)$ and numbers $(s,t)$, so that
$(D,s,t)\in\ParTheta^{-1}(\eta)$. It follows from our choice of
$\eta$, that the divisor is actually supported in
$\Sym^{g-1}(K)$. Moreover, looking on the $F$ side, we see
that the fiber points must converge to the divisor $q$. Thus, we see
that for all $T$ sufficiently large, the divisors in the fibers of
${\ParTheta'_T}^{-1}(\eta\times 0)$ are contained in the range of the
splicing map $$\Splice{T}\colon \Sym^{g-1}(\Sigma^c)\times F^c
\longrightarrow \Sym^{g}(\Sigma\#_T F),$$
where $\Sigma^c$ is a compact set whose interior contains $K$, and
$F^c$ is some compact subset of the punctured torus
(punctured at the connect sum point) which contains $q$.  

But applying Theorem~\ref{thm:Splicing}, we see that the maps 
$$\Sym^{g-1}(\Sigma^c)\times F^c \times [0,1]\times [0,1]
\longrightarrow \Targ\times \Jac(F)$$
obtained by mapping $$(D_1,D_2,s,t)\mapsto
\ParTheta'_T(\Splice{T}(D_1,D_2),s,t)$$ (which we will denote
$\ParTheta'_T\circ \Splice{T}$ in a mild abuse of notation) converge
in $C^1$ to the map which sends $$(D_1,D_2,s,t)\mapsto
\ParTheta(D_1,s,t)\times \AJ^{(1)}(D_2).$$ 
(Note that $\AJ^{(1)}(D_2)$ does not depend on $s$ and $t$ since we
are fixing the metric on the torus side.)  The preimage of $(\eta,0)$
under this limiting map is the fiber $$
(\ParTheta|_{\Sym^{g-1}(\Sigma^c)})^{-1}(\eta)\times q
=\ParTheta^{-1}(\eta)\times q$$ (the equality of the two sets follows
from the fact that $K\subset \Sigma_1^c$), which is used to calculate
$\theta$. Now the $C^1$ convergence, identifies this fiber with the
fiber $$(\ParTheta'_T\circ\Splice{T})^{-1}(\eta\times 0),$$ which is
used to calculate $\theta'(\spinc)$ (this is how we chose the subsets
$\Sigma_1^c\subset \Sigma_1$ and $F^c\subset F$). 
Thus, $\theta=\theta'$.

Note that our sign conventions are compatible with stabilization,
since if $\{\alpha_1,...,\alpha_g\}$, $\{\beta_1,...,\beta_g\}$ are
positively ordered for $U_0$ and $U_1$, then
$\{\alpha_1,...,\alpha_{g+1}\}$,
$\{\beta_{g+1},\beta_1,...,\beta_g\}$ are positively ordered for
$U_0'$ and $U_1'$, since the boundary of the two-cell corresponding to
$\beta_{g+1}$ is $-1$ times the boundary of the one-cell
corresponding to $\alpha_{g+1}$.
\end{proof}

We discuss another consequence of Theorem~\ref{thm:Splicing}, in a
case which will prove to be quite useful in the calculations.  But
first, we must characterize the image of the splicing map, in
terms of the connections. We content ourselves with a statement in the
case where $k_1=g_1-1$, as this is the only case we need to consider
in this paper.

\begin{prop}
\label{prop:RangeSplice}
Let $V_1\subset \Jac(\Sigma_1)$, $V_2\subset \Jac(\Sigma_2)$ be closed
subsets of the Jacobians. Suppose that $\Theta_{\met_1}^{-1}(V_1)$
contains no divisors which include the connect sum point $p_1$, and
suppose that $V_2$ contains no points in the theta divisor for
$\Sigma_2$.  Then, there are compact subsets $\Sigma_i^c\subset
\Sigma_i-p_i$ and a real number $T_0\geq 0$ so that for all $T\geq
T_0$, $\Theta_{h_T}^{-1}(V_1\times V_2)$ lies in the image of the
splicing map 
$$\Splice{T}\colon\Sym^{g_1-1}(\Sigma_1^c)\times\Sym^{g_2}(\Sigma_2^c)
\longrightarrow \Sym^{g-1}(\Sigma_1\#_T\Sigma_2).$$
\end{prop}

\begin{proof}
Our hypothesis on $V_1$ gives us a compact set
$K_1\subset \Sigma_1-p_1$ with the property that
$\Theta^{-1}_{\met_1}(V_1)\subset \Sym^{g-1}(K_1)$. Similarly, our
hypothesis on $V_2$ gives  a compact set $K_2\subset
\Sigma_2-p_2$ with the property that $(\AJ^{(g_2)})^{-1}(V_2)\subset
\Sym^{g_2}(K_2)$. We let $\Sigma_1^c$, $\Sigma_2^c$ be any pair of
compact sets whose interior contains $K_1$ and $K_2$.

Consider pairs $(A_T,\phi_T)$ over $\Sigma_1\#_T\Sigma_2$ which
correspond to the intersection of the theta-divisor with $V_1\times
V_2$, and which are normalized so that the $L^2$ norms over
$\Sigma_1\#_T\Sigma_2$ of $\phi_T$ is $1$.  By local compactness,
together with the fact that the tube admits no translationally
invariant harmonic spinor, any such sequence of pairs $(A_T,\phi_T)$
for tube-lengths $T\goesto\infty$ must admit a subsequence which
converges in $\Cinfty$ to an $L^2$ solution $\Phi_1$ and $\Phi_2$ on
the two sides $\Sigma_1^+$ and $\Sigma_2^+$, at least one of whose
$L^2$ norm is non-zero. By transfering back to $\Sigma_2$
(Proposition~\ref{prop:LTwoCohom}), our assumption on $V_2$ ensures
that $\Phi_2\equiv 0$.  By $\Cinfty$ convergence, then, the zeros of
$\phi_T$ must converge to the zeros of $\Phi_1$ over $\Sigma_1^+$.

Without loss of generality, we might as well assume that all the $A_T$
are of the form $A_1\#_T A_2$ for fixed $A_1\in\Jac(\Sigma_1)$,
$A_2\in\Jac(\Sigma_2)$. Note that for each $A_2\in V_2$, there is a
unique $A_2$-holomorphic section $\Phi_2$ over $\Sigma_2$ which, after
transferring to $\Sigma_2^+$, admits an asymptotic expansion 
$$\Phi_2 = e^{t/2} + \BigO(e^{-t/2})$$ 
(the growth here corresponds to the pole at $p_2\in \Sigma$ which we
have introduced in our convention for the Abel-Jacobi map).  Existence
of the section follows from the fact that the $g_2$-fold Abel-Jacobi
map has degree one (this is the ``Jacobi inversion theorem'', see
for instance p.~235 of~\cite{GriffHar}). Uniqueness
follows from the fact that a difference of two such would give an
$L^2$ section, showing that $A_2$ actually lies in the theta divisor,
which we assumed it could not.

We show the restrictions of $\phi_T$ to the $\Sigma_2$-side come close
to approximating $\Phi_2$ or, more precisely, that its zeros converge
to those of $\Phi_2$.

Rescale $\phi_T$ so that over $[-1,1]\times S^1$, it has the form
$$\phi_T = e^{-(T+t)/2} + \BigO(e^{-3(T+t)/2}).$$ Consider the section
$\Psi_T=\psi_2(\phi_T-e^{-T/2}\Phi_2)$, viewed as a section of
$\Sigma_2^+$ (we can do this, as its support is contained in the
support of $\psi_2$) . Note that $$\DBar_{A_2} \Psi_T = (\DBar \psi_2)
(\phi_T-e^{-T/2}\Phi_2).$$ Thus, $\|\DBar_{A_2}
\Psi_T\|=\BigO(e^{-3T/2})$. Since $\DBar_{A_2}$ is Fredholm with index
zero (it is a spin connection) and no kernel (it is not in the theta
divisor), it has no cokernel, and we can conclude that
$\|\Psi_T\|_{L^2(\Sigma_2^+)}\leq C e^{-3T/2}$ for some constant $C$
independent of $T$.

Since the restriction of $\Psi_T$ to $\Sigma_2^c$ is
$A_2$-holomorphic, elliptic regularity on this compact piece shows
that the section $\Psi_T =
\phi_T-e^{-T/2}\Phi_2$ is bounded by some quantity of order
$e^{-3T/2}$. Thus, the zeros of $\phi_T$ in $\Sigma_2^c$ converge to
those of $\Phi_2$.
\end{proof}

Armed with this proposition, we turn our attention to another
important consequence of Theorem~\ref{thm:Splicing}.  Let $\Sigma$ be
a surface of genus $g$, and let $\{\alpha_1,...,\alpha_g\}$ be a
complete set of attaching circles for a handlebody $U$ bounding
$\Sigma$. The holonomy around the first $g-1$ of the $\alpha_i$ gives
a map $$\Hol_{\alpha_1\times...\times\alpha_{g-1}}\colon
\Jac(\Sigma)\longrightarrow \Torus{g-1}.$$
According to~\cite{MacDonald} (see also~\cite{HigherType}, where a
related discussion is given), the preimage of a generic point in
$\Torus{g-1}$ via
$\Hol_{\alpha_1\times...\times\alpha_{g-1}}\circ\Theta_{\met}$ (for
any metric $\met$) is homologous to the torus
$\alpha_1\times...\times\alpha_{g-1}\subset \Sym^{g-1}(\Sigma)$. We
would like to find a metric on $\Sigma$ for which these spaces are
actually isotopic.

\begin{figure}
\mbox{\vbox{\epsfbox{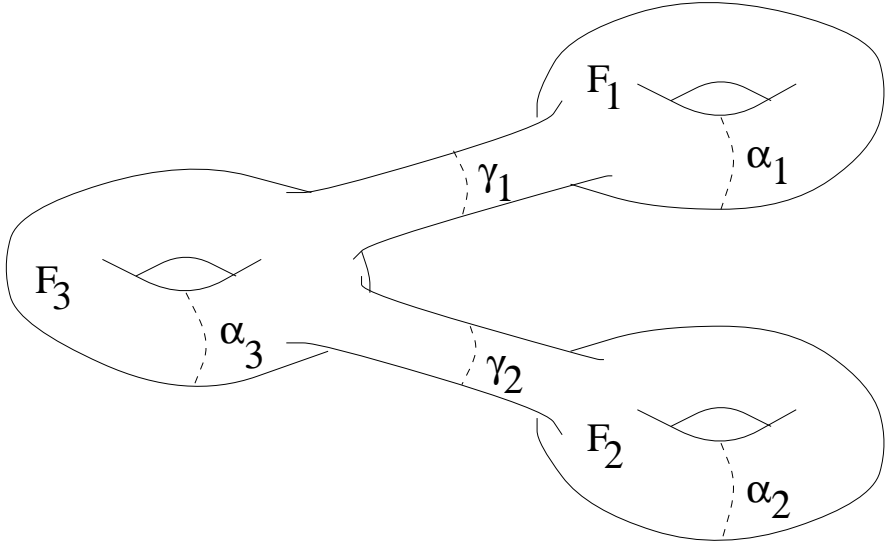}}}
\label{fig:Tori}
\caption{Connected sum of $F_1$ and $F_2$ with $F_3$. Attaching
circles $\{\alpha_1,\alpha_2,\alpha_3\}$ and connected sum circles
$\{\gamma_1, \gamma_2\}$ are included.} 
\end{figure}

To describe this metric, think of $\Sigma$ as a connected sum 
of $g-1$ disjoint tori $F_1,...,F_{g-1}$ with the remaining torus
$F_g$, in such a way that the curve $\alpha_i$ is supported in the
torus $F_i$ for $i=1,...,g-1$. Fix a metric $\met$ which is
product-like along the $g-1$ connect sum tubes, and let $\met(T)$
denote the metric obtained from $\met$ by stretching the connect sum
tubes by a factor of $T$. (The case where $g=3$ is illustrated in Figure~\ref{fig:Tori}.)

\begin{cor}
\label{cor:HolonomyConstraints}
Let $\Sigma$ be a surface of genus $g$ viewed as a connected sum of
tori as described above, and let $\{\alpha_1,...,\alpha_g\}$ be a
complete set of attaching circles. For any $\eta\in
\Torus{g-1}$  with the property that $\eta_i\neq 1/2$ for
all $i=1,...,g-1$, there is a $T_0$ so that for all $T\geq T_0$ the
subset
$(\Hol_{\alpha_1\times...\times\alpha_{g-1}}\circ\Theta_{\met(T)})^{-1}(\eta)$
is isotopic to the torus $\alpha_1\times...\times\alpha_{g-1}\subset
\Sym^{g-1}(\Sigma)$, where $\met(T)$ is the one-parameter family of
metrics obtained by stretching the connect sum tubes for the initial
metric $\met$.
\end{cor}

% old: % psfix -epsf -width 4  g3.eps > g3.ps
% new: psfix -epsf -width 5  g3.eps > g3.ps

\begin{proof}
We would like to apply a version of
Proposition~\ref{prop:RangeSplice}, with more than one neck (note that
the proof works in this context as well). Let $V_1$ be the theta
divisor of $F_g$. It contains none of the connect sum points, of
course, because it has degree zero. Moreover, the set
$\Hol_{\alpha_i}^{-1}(\eta_i)$ misses the theta divisor for $F_i$ for
$i=1,...,g-1$ (the theta divisor of $F_i$ consists of a single point
where the holonomy around $\alpha_i$ is $1/2$). Hence,
Proposition~\ref{prop:RangeSplice} applies: for all sufficiently long
necks, the theta divisor hits
$\Hol_{\alpha_1\times...\times\alpha_{g-1}}^{-1}(\eta)$ in a region
corresponding to the splicing map from
Theorem~\ref{thm:Splicing}. Thus, the composite of $\Theta$ with the
splicing map is $C^1$ close to the map $$F_1^c\times ... \times
F_{g-1}^c
\longrightarrow H^1(F_1; S^1)\times...\times
H^1(F_{g-1};S^1)\times H^1(F_g;S^1),$$ which is a product of $g-1$
copies of the Abel-Jacobi map with the inclusion of the point
(theta-divisor for the $F_g$). Since the Abel-Jacobi map in this case
is a diffeomorphism, the points where the $\alpha_i$-holonomy is
trivial forms a smoothly embedded curve. In fact, it is easy to see
that this curve is isotopic to $\alpha_i$ (see~\cite{MacDonald} and
also~\cite{HigherType}). Also, it is clear that post-composing with
evaluation along $\alpha_i$ gives us map to $(S^1)^{\times g-1}$ with
$\eta$ as a regular value, whose fiber is isotopic to
$\alpha_1\times...\times\alpha_{g-1}$. It is easy to see that any
other $C^1$ close map must have $\eta$ as a regular value, with an
isotopic fiber. Thus, the corollary follows from
Theorem~\ref{thm:Splicing}.
\end{proof}
\section{Calculations when $b_1(Y)>1$}
\label{sec:CalcBig}

The aim of this section is to prove the following:

\begin{theorem}
\label{thm:CalcBig}
When $b_1(Y)>1$, then the polynomial ${\underline\theta}$ is equal up
to sign to the symmetrized Alexander polynomial of $Y$.
\end{theorem}

In the proof of this theorem, we will naturally meet certain tori in
the symmetric product.  Given a Heegaard decomposition of $Y$, let
$\{\alpha_i\}$, $\{\beta_i\}$ be complete sets of attaching circles
for the two handlebodies. Given any $i$ and $j$, we have tori in
$\Sym^{g-1}(\Sigma)$ $$\Torus{}_i(\alpha)=\alpha_1\times..\times
{\widehat
\alpha}_i\times...\times\alpha_g$$
(where the notation indicates omission of the $i^{th}$ factor) and
$$\Torus{}_j(\beta)=\beta_1\times..\times {\widehat
\beta}_j\times...\times\beta_g.$$
We will show that the invariant $\tInv$ can be extracted from certain
polynomials associated to these tori; these polynomials are defined as
follows.  Let $\OurSym\rightarrow \Sym^{g-1}(\Sigma)$ be the covering
space of $\Sigma$ (as in Section~\ref{sec:DefTheta}) corresponding to
the kernel of the composite map:
$$\pi_1(\Sym^{g-1}(\Sigma))\rightarrow H_1(\Sym^{g-1}(\Sigma);\Z)\cong
H_1(\Sigma;\Z)\rightarrow H_1(Y;\Z)\cong H.$$ (Recall that $H$ is by
definition $H^2(Y;\Z)$.)  Thus, $H$ acts freely on $\OurSym$.  Let
${\AbTa{i}}$, ${\AbTb{j}}$ be a pair of lifts of $\Ta{i}$ and
$\Tb{j}$. Note that these lifts are tori, and indeed they map
isomorphically to $\Ta{i}$ and $\Tb{j}$ respectively. The intersection
points of $\Ta{i}$ with $\Tb{j}$ correspond to the intersection points
of ${\AbTa{i}}$ with the various translates under $H$ of the
torus $\AbTb{j}$. Then, we define a polynomial (an element of $\Z[H]$)
associated to the lifts $\AbTa{i}$ and $\AbTb{j}$ by the formula
\begin{equation}
\label{eq:DefOfC}
C_{i,j}=\sum_{h\in H}
\#\left(\AbTa{i}\cap h\AbTb{j}\right) [h].
\end{equation}
(For the intersection numbers here, we use orientations for the
$\AbTa{i}$ and $\AbTb{j}$ induced from orientations of $\Ta{i}$ and
$\Tb{j}$; we return to a more careful discussion of signs in
Section~\ref{subsec:Refinements}.) Summing over the action of $\Tors$,
we get an induced polynomial ${\underline C}_{i,j}\in
\Z[\tH]$ (recall that $\tH=H/\Tors$). In Section~\ref{sec:Alex}, we
will show that the Alexander polynomial of $Y$ is the greatest common
divisor of the ${\underline C}_{i,j}$ for $i,j=1,...,g$. Different
lifts of the $\Ta{i}$ and $\Tb{j}$ give rise to translates of the
$C_{i,j}$ and ${\underline C}_{i,j}$ by elements in $H$.

The main ingredient in the proof of Theorem~\ref{thm:CalcBig} is a
perturbation of the invariants, which corresponds to moving the tori
$L_0$ and $L_1$. Let $\La{i}$ be the space of $B\in\Jac(\Sigma)$ with
$\Hol_{\alpha_k}B=0$ for all $k\neq i$, and similarly let $\Lb{j}$ be
the space of $B\in\Jac(\Sigma)$ with $\Hol_{\beta_k}B=0$ for all
$k\neq j$. We will move the tori $L_0$ and $L_1$ inside $\La{i}$ and
$\Lb{j}$, and obtain an expression of $\theta$ in terms of the 
intersection of the tori $\La{i}$ and $\Lb{j}$ with the theta divisor. 
An important point, then, is that we can concretely
understand these intersections, 
for favorable initial metrics. It is with the help of this
description, then, that we meet the polynomials described above.
But first, we describe how to calculate $\theta$ in terms of the lifts of
$\La{i}$ and $\Lb{j}$. To do this, we discuss the lifts in detail.

There is a lift ${\widetilde\delta}\colon \OurJac\longrightarrow
H^2(Y;\R)$ of the coboundary map $H^1(\Sigma;S^1)\rightarrow
H^2(Y;S^1)$, which is uniquely specified once we ask that
${\widetilde\delta}(\Lag_0(\spinc))=0$. With our conventions, then,
${\widetilde\delta}(\Lag_1(\spinc))={\underline\spinc}$. 

To $\alpha_i$, assign an element $\alpha_i^*\in H^1(\Sigma;\Z)$, as
follows. Let $\gamma_i$ be the core of the $i^{th}$ one-handle in
$U_0$ (i.e. this is the oriented curve which intersects only the
attaching disk associated to $\alpha_i$, which it intersects
positively in a single point), then $\alpha_i^*$ is the Poincar\'e
dual (in $\Sigma$) to a class in $H_1(\Sigma;\Z)$ whose image in
$H_1(U;\Z)$ is represented by $\gamma_i$.  (The class $\alpha_i^*$ is
not uniquely determined by this property, but the our constructions
involving $\alpha_i^*$ are independent of its choice.) Note that the
class $\mu_i=\delta\alpha_i^*\in H^2(Y;\Z)$, is Poincar\'e dual (in
$Y$) to the homology class represented by $\gamma_i$. The element
$\beta_i^*$ is defined in the analogous manner, only using $U_1$
instead of $U_0$. We let $\nu_j$ denote $\delta\beta_j^*$.  Choose $i$
and $j$ so that $\mu_i$ and $\nu_j$ are not torsion classes; we can
find such $i$ and $j$ since $H_1(U_0)$ and $H_1(U_1)$ both surject
onto $H_1(Y)$.

By multiplying $\alpha_i$ by $-1$ if necessary, we can assume that
${\underline\mu}_i$ and ${\underline\nu}_j$ in $H^2(Y;\Z)/\Tors$ are
not negative multiples of each other. (In keeping with the conventions
introduced in Section~\ref{sec:Introduction}, we underline objects
when viewing them modulo the action of torsion.) We define subsets
$\Lambda_0^+(\spinc),\Lambda_1^+(\spinc)\subset \OurJac$ which
correspond to all translates of (small perturbations of) $\Lag_0$ and
$\Lag_1$ in the directions determined by $\alpha_i^*$ and $\beta_j^*$;
more precisely, 
\begin{eqnarray*}
\Lambda_0^+(\spinc)=\Lag_0(\spinc)+\eta_0+ \R^+{\alpha_i^*}
&{\text{and}}&
\Lambda_1^+(\spinc)=\Lag_1(\spinc)+\eta_1- \R^+{\beta_j^*}.
\end{eqnarray*}
Under the map $\OurJac\rightarrow \Jac$, the spaces
$\Lambda_0^+(\spinc)$ and $\Lambda_1^+(\spinc)$ project to $\La{i}$
and $\Lb{j}$ respectively.

In the case where $b_1(Y)=2$, we make use of special allowable metrics:

\begin{defn}
Fix an $1/2 > \epsilon>0$, and let $U$ be a handlebody which bounds
$\Sigma$.  A metric $\met$ on $\Sigma$ is said to be {\em strongly
allowable for $\epsilon$} if it is product-like in a neighborhood of
$g$ attaching circles $\{\gamma_1,...,\gamma_g\}$, and any point in
the theta divisor for $\Sigma$ must have holonomy around some
attaching circle $\gamma_i$ within $\OneHalf\epsilon$ of $\OneHalf$.
\end{defn}

Given any $\epsilon>0$, there exist metrics which are 
strongly allowable for $\epsilon$ thanks to Lemma~\ref{lemma:MissTheta}.

\begin{prop}
\label{prop:Intersection}
Let $\spinc$ be any $\SpinC$ structure on $Y$, and fix
$\Lambda^+_i(\spinc)$ for $i=0,1$ as above -- using classes
${\underline\mu}_i$ and ${\underline\nu}_j$ which are not negative
multiples of one another. There is an $\epsilon>0$ with the property
that for any metrics $\met_0$ and $\met_1$ which are $U_0$ and
$U_1$-allowable respectively, where $\met_0$ is also
$\epsilon$-strongly $U_0$-allowable then, we can find $\eta_0$,
$\eta_1$ sufficently small, with $$\theta(\spinc)=
\#\left(\OurTheta^{-1}_{\met_0}\left(\Lambda_0^+(\spinc)+\eta_0\right)
	\cap \OurTheta^{-1}_{\met_1}\left(\Lambda_1^+(\spinc)+\eta_1\right)\right).
$$
\end{prop}

\begin{proof}
The proof will rely on the fact that $\OurDelta\circ \OurTheta_{\met_t}$
has bounded variation along any one-parameter family of
metrics. Specifically, let $\met_t$ be a one-parameter family of
metrics, and fix a norm on $H^2(Y;\R)$. Then, there is a constant $K$
with the property that for any $D\in\OurSym$, $s,t\in [0,1]$,
\begin{equation}
\label{eq:BoundedVariation}
|\OurDelta\circ\OurTheta_{\met_s}(D)-\OurDelta\circ\OurTheta_{\met_t}(D)|<
K.
\end{equation} This follows immediately from the compactness of
$\Sym^{g-1}(\Sigma)$, together with the fact that
$\OurDelta\circ\OurTheta$ is an $H^2(Y;\Z)$-equivariant map.

Note that $\theta(\spinc)$ is calculated by the number of points
(counted with signs) in the zero-dimensional submanifold of
$\OurSym\times[0,1]\times[0,1]$ $$\{(D,s,t)\big| s\leq t,
\OurTheta_{\met_s}(D)\in
\Lag_0(s)+\eta_0, \OurTheta_{\met_t}(D)\in\Lag_1(\spinc)+\eta_1\},$$ 
a space we denote by $M(\spinc)$. We construct a cobordism between
this space and the points in the intersection stated in the lemma, as
follows.  We can assume without loss of generality that $\met_t$ is
constant between $[0,1/4]$ and $[3/4,1]$. Moreover, let $\psi$ be a
non-decreasing smooth function on $[0,1]$ which is monotone increasing
in the range $[0,1/4]$, with $\psi(0)=0$ and $\psi|_{[1/4,1]}\equiv
1$. Consider the subspace of $\OurSym\times[0,1]\times[0,1]$ (which
agrees with $M(\spinc)$ when $u_1=u_2=0$): $$M_{u_1,u_2}(\spinc)=
\left\{
(D,s,t)\Bigg|
\begin{array}l
s\leq t, \\
\OurTheta_{\met_s}(D)\in
\Lag_0(\spinc)+\eta_0+u_1\psi(s) \alpha_i^*, \\
\OurTheta_{\met_t}(D)\in\Lag_1(\spinc)+\eta_1-u_2\psi(1-t) \beta_j^*
\end{array}
\right\}.$$

We argue that for all sufficiently large $u$,
\begin{equation}
\label{eq:PerturbedPoints}
M_{u,u}=\OurTheta^{-1}_{\met_0}\left(\Lambda_0^+(\spinc)+\eta_0\right)\cap
\OurTheta^{-1}_{\met_1}\left(\Lambda_1^+(\spinc)+\eta_1\right).
\end{equation}
Since ${\underline\mu}_i$ and ${\underline\nu}_j$ are not negative
multiples of one another, we see that that as $u\goesto\infty$, the
distance between the point ${\widetilde\delta}(L_0(\spinc)+\eta_0+u
\alpha_i^*)\in H^2(Y;\R)$ and the ray
${\widetilde\delta}(L_1(\spinc)+\eta_1-\R^+\beta_j^*)$ goes to
infinity, and similarly the distance between the point
${\widetilde\delta}(L_1(\spinc)+\eta_1-u\beta_j)$ and the ray
${\widetilde\delta}(L_0(\spinc)+\eta_0+\R^+\alpha_i^*)$ goes to
infinity. Fix $u$ large enough that both distances are larger than the
constant $K$ from Inequality~\eqref{eq:BoundedVariation}. This
condition ensures that all points $(D,s,t)\in M_{u,u}$ have $s\leq
1/4$ and $t\geq 3/4$. Monotonicity of $\psi$ over $[0,1/4]$, and the
choice of $u$ then also ensures that the identification
\eqref{eq:PerturbedPoints} holds.

Thus, Proposition~\ref{prop:Intersection} is established once we
construct a smooth cobordism between $M_{0,0}$ and $M_{u,u}$.
Consider the spaces obtained by connecting $M_{0,0}$ to $M_{u,u}$ by
first allowing $u_1$ to go from $0$ to $u$ (to connect $M_{0,0}$ to
$M_{u,0}$) and then allowing $u_2$ to go from $0$ to $u$ (to connect
$M_{u,0}$ to $M_{u,u}$). Since $\met_0$ and $\met_1$ are allowable
metrics, the $s=0$ and $t=1$ boundaries are excluded in this
one-parameter family for all small $\eta_0$ and $\eta_1$.  Thus, we
get a cobordism between $M_{0,0}$ and $M_{u,u}$, provided that the
$M_{u_1,u_2}$ do not hit the $s=t$ boundary, which is guaranteed if
$$\big(\La{i}+\eta_0\big)\cap
\big(\Lb{j}+\eta_1\big)=\emptyset. $$ Taking $\OurDelta$ of both
spaces, we get a pair of lines in $H^1(Y;\R)$, which generically miss
each other when $b_1(Y)>2$.

The case where $b_1(Y)=2$ requires a slightly closer investigation.
In the first part of the cobordism, where we allow $u_1$ to vary in
$M_{u_1,0}$, there are still no $s=t$ boundaries, as we can arrange
for $$(\La{i}+\eta_0)\cap (\Lag_1+\eta_1)=\emptyset$$ (since, applying
$\OurDelta$, we have a point and a line in a two-space).  Now, as
$u_2$ varies in the $M_{u,u_2}$, it is easy to see that the only
possible $s=t$ boundaries lie in the range where $s\leq 1/4$, by our
hypothesis on $\psi$, and hence they must lie in the set
$$\Theta_{\met_0}^{-1}\big((\La{i}+\eta_0)\cap (\Lb{j}+\eta_1)\big),$$
since $\met_t$ is constant for $t\leq 1/4$.  Now, consider the
intersection point $p$ of the induced rays $\OurDelta(\La{i}+\eta_0)$
and $\OurDelta(\Lb{j}+\eta_1)$.  Note that as $\met_0$ is stretched
out normal to the attaching disks, the image under $\OurDelta$
of the intersection
intersection $$(\La{i}+t\alpha_i^*)\cap
\Theta_{\met_0}(\Sym^{g-1}(\Sigma))$$
converges to a discrete subset of the ray
$\R^+{\underline\mu}_i\subset H^2(Y;\R)$, consisting of points
separated by some distance $\delta>0$ (which depends on the $\mu_i$
and $\nu_j$). If $p$ misses this discrete set, then if $\met_0$ is
sufficiently stretched out, then all sufficiently small $\eta_0$ and
$\eta_1$ have the property that
$$\Theta_{\met_0}^{-1}\big((\La{i}+\eta_0)\cap
(\Lb{j}+\eta_1)\big)=\emptyset.$$
If, on the other hand, $p$ lies on the discrete set, then, given any
sufficiently small
$0<\gamma$, if $\met_0$ is
sufficiently stretched out, then all sufficiently small $\eta_0$ and
$\eta_1$ have the property that
$$\Theta_{\met_0}^{-1}\big((\La{i}+\eta_0+\gamma \alpha_i^*)\cap
(\Lb{j}+\eta_1)\big)=\emptyset.$$
Thus, in both cases, we have obtained the requisite cobordism.
\end{proof}

Moreover, we can describe the intersection appearing in
Proposition~\ref{prop:Intersection}, in terms of the tori $\Ta{i}$ and
$\Tb{j}$ described in the beginning of the section.  But first, we
state a relevant lemma, whose proof fits naturally into the framework
of Section~\ref{sec:MoveTheta}.

\begin{lemma}
\label{lemma:StretchEverything}
Let $\Sigma$ be a surface, realized as a connected sum of $g$ tori as
in Corollary~\ref{cor:HolonomyConstraints}, and let
$\{\alpha_1,...,\alpha_g\}$ a complete set of attaching circles for
the handlebody $U$ which bounds $\Sigma$, each of which is disjoint
from the separating curves for the connected sum decomposition of
$\Sigma$. Fix a metric $h$ which is product-like near the $\alpha_i$
and the separating curves. Then, there is a $T_0$ so that any metric
which is obtained from $h$ by stretching each of the $2g-1$ curves at
least by $T_0$ is $U$-allowable.
\end{lemma}

\begin{proof}
Take a weak limit of connections which lie in $L(U)$, as all the
$2g-1$ curves are stretched. Under this limit, the surface degenerates
into a collection of genus zero surfaces (with cylindrical ends),
whose ends correspond to attaching circles for $U$ or separating
curves for $\Sigma$. Thus, the weak limit of connections in $L(U)$
induces a connection over these genus zero surfaces, whose holonomies
around all its bounding circles is zero.  But none of these support
harmonic spinors according to Proposition~\ref{prop:LTwoCohom},
proving the lemma.
\end{proof}

\begin{prop}
\label{prop:IdentTori}
There is a $U_0$-allowable metric $h$, for which
$\Theta_{h}^{-1}(\La{i}+\eta)$ is isotopic to $\Ta{i}$ for all
generic, small $\eta$.
\end{prop}

\begin{proof}
Fix a metric $\met$ on $\Sigma$ as in
Lemma~\ref{lemma:StretchEverything}.  Let $\met(T)$ denote the
metric which is stretched by $T_0$ along the $\alpha_i$ and $T$ along
the separating curves. The lemma guarantees
that for all $T>T_0$, $\met(T)$ is allowable. Now, for all
sufficiently large $T$, Corollary~\ref{cor:HolonomyConstraints} gives
the isotopy of $\Ta{i}$ with the subset $\Theta_{h}^{-1}(\La{i})$,
where $h=\met(T)$. 
\end{proof}

Putting together Propositions~\ref{prop:Intersection} and
\ref{prop:IdentTori}, we obtain the following topological
description of $\theta$.

\begin{prop}
\label{prop:CalcTheta}
For some $\SpinC$ structure $\spinc\in\SpinC(Y)$, we have that
$$C_{i,j}[\spinc]=(1-\mu_i)(1-\nu_j)\cm\theta.$$
\end{prop}

\begin{proof}
Let $\met_k$ be $U_k$-allowable metrics for $k=0,1$.  Given $i,j$, we
construct a natural element ${\widetilde C}_{i,j}\in\Z[\SpinC(Y)]$
closely related to the $C_{i,j}$ defined in the beginning of this
section. The ${\widetilde C}_{i,j}$ will be a translate of the
following analogue of the $C_{i,j}$, which is assigned to a pair of
lifts ${\widetilde A}_i$ and ${\widetilde B}_j$ (to $\OurSym$) of the
manifolds $A_i=\Theta^{-1}_{\met_0}(\La{i}+\eta_0)$ and
$B_j=\Theta^{-1}_{\met_1}(\Lb{j}+\eta_1)$: $$ X_{i,j}= \sum_{h\in H}
\# \left({\widetilde B}_j \cap h {\widetilde A}_i\right) [h].
$$ To do define the ${\widetilde C}_{i,j}$, we must assign a $\SpinC$
structure to each intersection point of $A_i$ with $B_j$. 
To this end,
we assume that $\met_0$ and $\met_1$ are strongly allowable for some
$\epsilon>0$.  If $p\in A_i\cap B_j$, then let ${\widetilde p}$ be a
lift of $p$. There is a pair of lifts ${\widetilde A}_i$ and
${\widetilde B}_j$ of $A_i$ and $B_j$ which meet in ${\widetilde
p}$. 
Note that there are unique lifts $\Lag_0$ and $\Lag_1$
whose image under $\OurDelta$ lie
in an $\epsilon$ neighborhood of ${\widetilde p}-\OneHalf{\underline\mu}_i$
and ${\widetilde p}+\OneHalf{\underline\nu}_j$ respectively.
Let $G(p)$ denote the $\SpinC$ structure which corresponds
to the difference between these two lifts (i.e. if
${\widetilde\Lag}_0$ and ${\widetilde\Lag_1}$ are the two lifts, then
the pair $L_0(G(p))$ and $L_1(G(p))$ are translates of
${\widetilde\Lag}_0$ and ${\widetilde\Lag_1}$ by a single cohomology
class in $H^2(Y;\Z)$).  By summing over all intersection points which
correspond to a given $\SpinC$ structure (with signs), we obtain the
element ${\widetilde C}_{i,j}\in\Z[\SpinC(Y)]$, which is clearly the
translate by some $\SpinC$ structure of the polynomial $X_{i,j}\in
\Z[H]$ defined above.

It follows from Proposition~\ref{prop:Intersection} that
$\theta(\spinc)$ is given by adding up the intersection number of
certain lifts of $A_i$ and
$B_j$. Moreover, the intersection point $p$ will
contribute for each $k,\ell\geq 0$ in the $\SpinC$ structure
$$G(p)+k\mu_i+\ell\nu_j$$ (as those are the $\SpinC$ structures for
which $p$ lies on the corresponding rays).  This proves that
\begin{equation}
\label{eq:CalculateTheta}
\theta={\widetilde C}_{i,j}\left(\sum_{k=0}^{\infty}\mu_i^k\right)
\left(\sum_{\ell=0}^{\infty}\nu_j^\ell\right).
\end{equation}

Finally, the proposition is proved once we establish that the
${\widetilde C}_{i,j}$ is a translate of $C_{i,j}$, as defined in
Equation~\eqref{eq:DefOfC}.  To see this, recall that
Proposition~\ref{prop:IdentTori} guarantees that the spaces $A_i$ and
$B_j$ are isotopic, for suitable choices of allowable metrics $\met_0$
and $\met_1$, to $\Ta{i}$ and $\Tb{j}$ respectively.  Now the
polynomial $X_{i,j}$, which is clearly a translate of ${\widetilde
C}_{i,j}$, depends on the submanifolds $A_i$ and $B_j$ only up to
isotopy.  Thus, the proposition follows.
\end{proof}

In particular, we have the following:

\begin{cor}
If ${\underline\mu}_i$ and ${\underline\nu}_j$ are not negative
multiples of one another, then $${\underline
C}_{i,j}=(1-{\underline\mu}_i)(1-{\underline\nu}_j){\underline\theta}.$$
\end{cor}

Theorem~\ref{thm:CalcBig} follows from this corollary. First, after
handle-slides, we can arrange that all the ${\underline\mu}_i$ and
${\underline\nu}_j$ are non-zero, so that ${\underline\theta}$ divides
each ${\underline C}_{i,j}$ and hence the Alexander polynomial.
Furthermore, after additional handleslides we can arrange that
${\underline\mu}_1={\underline\nu}_1$,
${\underline\mu}_2={\underline\nu}_2$, and ${\underline\mu}_1$ and
${\underline\mu}_2$ are linearly independent in $H^2(Y;\R)$. This
shows that the greatest common divisor of ${\underline C}_{1,1}$ and
${\underline C}_{2,2}$ is ${\underline\theta}$, so the latter agrees
with the Alexander polynomial.

\section{Calculating the invariant when $b_1(Y)=1$}
\label{sec:CalcOne}

The aim of this section is to prove the following:

\begin{theorem}
\label{thm:CalcOne}
Let $A=a_0+\sum_{i=1}^k a_i(T^i+T^{-i})$ be the symmetrized Alexander
polynomial of $Y$, normalized so that $A(1)=|\Tors H_1(Y;\Z)|$. Then,
the Laurent polynomial of ${\underline\theta}$ is equal to
$${\underline\theta}=b_0+\sum_{i=1}^{\infty} b_i (T^i+T^{-i}),$$ where
$$b_i=\sum_{j=1}^\infty j\cdot a_{i+j}.$$
\end{theorem}

We use the same notation as in Section~\ref{sec:CalcBig}. Again, we
need that ${\underline\mu}_i$ and ${\underline\nu}_j$ are non-zero. A
bit more care is needed in defining the $\Lambda_i^+(\spinc)$. By
multiplying the $\mu_i$ and $\nu_j$ by $-1$ if necessary, we can
arrange that ${\underline\mu}_i$ and ${\underline\nu}_j$ are positive
multiples of each other. Indeed, to simplify the language, we can
choose an isomorphism $H^2(Y;\R)\cong \R$ so that ${\underline\mu}_i$
and ${\underline\nu}_j$. Suppose that, under this identification,
${\underline\spinc}\leq 0$. In this case, we define $\Lambda_i^+$ for
$i=0,1$ as in Section~\ref{sec:CalcBig}. The proof of
Proposition~\ref{prop:Intersection} applies to give us the analogous
result:

\begin{prop}
\label{prop:IntersectionOne}
Fix an identification $H^2(Y;\R)\cong \R$ and suppose that with
respect to this identification, the classes ${\underline\mu}_i$ and
${\underline\mu}_j$ are both positive. Let ${\spinc}$ be a metric so
that ${\underline\spinc}\leq 0$. Then, if
$\met_0$ and $\met_1$  are $U_0$- resp. $U_1$-allowable metrics, 
then we have the relation:
$$\theta(\spinc)=
\#\left(\OurTheta^{-1}_{\met_0}\left(\Lambda_0^+(\spinc)+\eta_0\right)
	\cap
\OurTheta^{-1}_{\met_1}\left(\Lambda_1^+(\spinc)+\eta_1\right)\right)  $$
for all sufficiently small, generic $\eta_0$ and $\eta_1$. 
\end{prop}

\begin{proof}
We adopt the notation and most of the argument for the proof of
Proposition~\ref{prop:Intersection}. The difference arises when one
wishes to prove that the moduli spaces $M_{0,0}$ and $M_{u,u}$ are
cobordant, i.e. when one wishes to exclude the possible $s=t$ boundary
components. To do this, it is no longer possible to argue that
$$(\Lambda_i(\alpha)+\eta_0)\cap (\Lag_1+\eta_1)=\emptyset.$$ Rather,
to exclude $s=t$ boundaries, we show that
$$\left(\Lag_0(\spinc)+\eta_0+u_1\psi(s)\alpha_i^*\right)\cap
\left(\Lag_1(\spinc)+\eta_1-u_2\psi(1-s)\beta_j^*\right)=\emptyset$$
(which suffices). To see this, begin by choosing generic $\eta_0$ and
$\eta_1$ so that $\OurDelta(\eta_0-\eta_1)$ is positive (we are free
to do this according to Proposition~\ref{prop:IndepOfMetOne}, or just
Proposition~\ref{prop:IndepOfMetOneZ} for integral homology
three-spheres).  If the intersection were non-empty, by applying
$\OurDelta$ to both sides, we would get: $$\OurDelta(\eta_0) + u_1
{\underline\mu}_i = {\underline \spinc} +
\OurDelta(\eta_1)-u_2 {\underline\nu_j}.$$ 
Thus, $$\OurDelta(\eta_0-\eta_1) + u_1 {\underline\mu}_i - {\underline
\spinc} + u_2 {\underline\nu}_j=0,$$ which is impossible, as it is a sum
of four non-negative terms at least one of which (the first) is
positive.
\end{proof}

\begin{remark}
The hypothesis that ${\underline\spinc}$ has an opposite sign from
${\underline\mu}_i$ and ${\underline\nu}_j$ is important. If it is
violated, there will be addition correction terms from $s=t$ boundary
components in the cobordism.
\end{remark}

We define ${\widetilde C}_{i,j}$ as before.

\begin{prop}
\label{prop:CalcPartOne}
For all $\spinc$ for which ${\underline\spinc}$ is a non-negative
multiple of ${\underline\mu_i}$, the value of $\theta(\spinc)$ is
the coefficient of $[\spinc]$ in the Laurent series
$${\widetilde C}_{i,j}\left(\sum_{k=0}^{\infty}[\mu_i]^k\right)
\left(\sum_{\ell=0}^{\infty}[\nu_j]^\ell\right).$$
\end{prop}

To finish the proof of Theorem~\ref{thm:CalcOne}, after handleslides,
we arrange that ${\underline\alpha}_g^*={\underline\beta}_g^*$ is a
generator of $\tH$. Let $\tC$ be the image of 
${\widetilde C}_{g,g}[-\alpha_g^*]$ in $\Z[\tH]$. Write
$$\tC=\sum_{k\in\Z}d_k T^k,$$ where $T$ corresponds to the generator of
$\Z[\tH]$. Then:

\begin{prop}
\label{prop:CalcPart}
Write $${\underline \theta}=b_0 + \sum b_i (T^i+T^{-i}).$$ Then,
$$b_i=\sum_{j=1}^\infty j d_{i-j},$$ and also 
$$b_i=\sum_{j=1}^\infty j d_{j-i}.$$ 
\end{prop}

\begin{proof}
This is a direct consequence of Proposition~\ref{prop:CalcPartOne}.
\end{proof}

\begin{prop}
$\tC$ is the symmetrized Alexander polynomial of $Y$.
\end{prop}

\begin{proof}
Proposition~\ref{prop:CalcPart} shows that $\tC$ is determined by
$\tInv$ and the classes ${\underline\mu}_i$ and ${\underline\mu}_j$.
After a series of handleslides, we can arrange that all
${\underline\mu}_i$ and ${\underline\nu}_j$ are equal to one, fixed
generator of $\tH$. It follows from Proposition~\ref{prop:CalcPart}
that all ${\underline C}_{i,j}$ equal $\tC$ up to translation.  Since
the Alexander polynomial $A$ (modulo multiplication by $\pm T^i$) is
the greatest common divisor of the ${\underline C}_{i,j}$
(c.f. Proposition~\ref{prop:Alex}), it follows that $\tC$ is a
translate of the symmetrized Alexander polynomial. But
Proposition~\ref{prop:CalcPart} also shows that $\tC$ is symmetric.
\end{proof}

\section{The Alexander Polynomial}
\label{sec:Alex}

In the calculation of the invariant, we have met certain tori in the
symmetric product, to which we associated polynomials $C_{i,j}$ and
$\tC_{i,j}$.  The aim of this section is to relate them to the
Alexander polynomial and, in Subsection~\ref{subsec:Refinements}, to
relate them to Turaev's torsion invariant.

Recall that a Heegaard decomposition of $Y$ and complete sets of attaching circles $\{\alpha_i\}$, $\{\beta_i\}$
for the two handlebodies naturally give rise to tori,
indexed by $i,j=1,...,g$ 
$$\Torus{}_i(\alpha)=\alpha_1\times..\times {\widehat
\alpha}_i\times...\times\alpha_g$$
and 
$$\Torus{}_j(\beta)=\beta_1\times..\times {\widehat
\beta}_j\times...\times\beta_g$$
in $\Sym^{g-1}(\Sigma)$.
Let $\OurSymF\rightarrow \Sym^{g-1}(\Sigma)$ be the covering space of
$\Sigma$ corresponding to the kernel of the composite map:
$$\pi_1(\Sym^{g-1}(\Sigma))\longrightarrow H_1(\Sym^{g-1}(\Sigma))\cong
H_1(\Sigma)\rightarrow H_1(Y)/\Tors=\tH.$$ Thus, $\tH$ acts freely on
$\OurSymF$. 
Let ${\AbTa{i}}$, ${\AbTb{j}}$ be a pair of lifts of
$\Ta{i}$ and $\Tb{j}$. Note that these lifts are tori, and indeed they
map isomorphically to $\Ta{i}$ and $\Tb{j}$ respectively. The
intersection points of $\Ta{i}$ with $\Tb{j}$ correspond to the
intersection points of ${\AbTa{i}}$ with the various translates under
$\TorsFree$ of the torus $\AbTb{j}$. Then, we define a polynomial (an
element of $\Z[\TorsFree]$) associated to the lifts $\AbTa{i}$ and
$\AbTb{j}$ by the formula 
$$\tC_{i,j}=\sum_{h\in \TorsFree}
\#\left(\AbTb{j}\cap h\AbTa{i}\right) [h].$$
Note that this agrees with the earlier definition of ${\underline
C}_{i,j}$.

\begin{prop} 
\label{prop:Alex}
The Alexander polynomial of $Y$ is the greatest common
divisor of the ${\underline C}_{i,j}$ for $i=1,...,g$.
\end{prop}

Before giving the proof, we briefly recall how to calculate the
Alexander polynomial from a Heegaard decomposition. The Heegaard
decomposition gives rise to a CW complex structure on $Y$, with one
zero-cell in $U_0$, $g$ one-cells (one for each handle in $U_0$;
i.e. these are obtained by pushing $g$ curves over $\Sigma$ which are
dual to the attaching circles $\{\alpha_1,...,\alpha_g\}$), $g$
two-cells (attached to $\Sigma$ along the attaching circles
$\{\beta_1,...,\beta_g\}$), and one three-cell. Let ${\widetilde Y}$
be the maximal free Abelian cover of $Y$, i.e. the one corresponding
to the kernel of $$\pi_1(Y)\rightarrow \tH.$$ This space
inherits a natural action of $\TorsFree=H_1(Y;\Z)/\Tors$.  Moreover,
the lifts of the cells in $Y$ gives and $\TorsFree$-equivariant
CW-complex structure on ${}$ ${\widetilde Y}$; more precisely, choose
for each cell in $Y$ a single cell in ${\widetilde Y}$ which covers it
(this is a {\em fundamental family of cells} in the sense
of~\cite{Turaev}). Then these cells form a basis of the chain complex
$C_*({\widetilde Y})$ over the group-ring $\Z[\TorsFree]$.  Thus, we
can view the cellular boundary from two-chains to one-chains on
${\widetilde Y}$ as a map $$\partial\colon
(\Z[\TorsFree])^g\longrightarrow (\Z[\TorsFree])^g;$$ i.e. it is a
$g\times g$ matrix over $\Z[\TorsFree]$. Given $i,j$, let
$\tDelta_{i,j}$ be the determinant of the $(g-1)\times (g-1)$ minor
obtained by deleting the $i^{th}$ row and the $j^{th}$ column of this
matrix. The {\em Alexander polynomial of $Y$}, then, is the greatest
common divisor of the $\tDelta_{i,j}$. 

\vskip0.3cm
\noindent
{\bf{Proof of Proposition~\ref{prop:Alex}.}}
Let ${\widetilde \Sigma}$ denote the cover of $\Sigma$ corresponding
to the kernel of the map $\pi_1(\Sigma)\rightarrow \tH$. Note that
the space $\OurSymF$ is the 
quotient of $\Sym^{g-1}({\widetilde\Sigma})$ by the equivalence
relation
$$\{x_1,...,x_{g-1}\}\sim
\{h_1 x_1,...,h_{g-1} x_{g-1}\}$$
for all tuples $(h_1,...,h_{g-1})\in \TorsFree^{g-1}$ with
$\sum_{i=1}^{g-1}h_i=0$.

Now, let $\{a_1,...,a_g\}$ be the one-cells corresponding to the
$\{\alpha_1,...,\alpha_g\}$, and $\{b_1,...,b_g\}$ be the two-cells
corresponding to $\{\beta_1,...,\beta_g\}$. Let $\{{\widetilde
\alpha}_1,...,{\widetilde
\alpha}_g\}$ and  $\{{\widetilde \beta}_1,...,{\widetilde \beta}_g\}$
be lifts of the corresponding attaching circles in ${\widetilde
\Sigma}$; and $\{{\widetilde a}_1,...,{\widetilde a}_g\}$ and
$\{{\widetilde b}_1,...,{\widetilde b}_g\}$ denote the corresponding
lifts of the associated cells in ${\widetilde Y}$.  Then, the formula
for the boundary map is given by $$
\partial {\widetilde
b}_i = \sum_{j=1}^g \left(\sum_{h\in \TorsFree} \#\left({\widetilde
\beta}_i\cap h{\widetilde \alpha}_j\right)[h]\right) {\widetilde
a}_j.$$

From this, then, we can obtain the identification of ${\underline
C}_{i,j}=\tDelta_{i,j}$. For notational convenience, we write this out
for $i=j=g$, but the general case follows in the same manner:
\begin{eqnarray*}
{\underline C}_{g,g} &=&
(-1)^{g-1}
\sum_{h\in \TorsFree} \#\left(\AbTb{g}\cap h\AbTa{g}\right) [h] \\
&=&
(-1)^{g-1}
\sum_{h\in \TorsFree}  \left(
\sum_{h_1+...+h_{g-1}=h}
\#\left((\AbBeta_1\times...\times\AbBeta_{g-1})\cap 
(h_1 \AbAlpha_1\times...h_{g-1}\AbAlpha_{g-1})\right)\right)[h] \\
&=& 
(-1)^{g-1}
\sum_{h_1,...,h_{g-1}}
\#\left((\AbBeta_1\times...\times\AbBeta_{g-1})\cap 
(h_1 \AbAlpha_1\times...h_{g-1}\AbAlpha_{g-1})\right)[h_1+...+h_{g-1}] 
\end{eqnarray*}
\begin{eqnarray*}
&=& (-1)^{\epsilon}\sum_{h_1,...,h_{g-1}}\sum_{\sigma\in \Perm{g-1}} (-1)^{\sigma}
\#\left(\AbBeta_1\cap
h_1\AbAlpha_{\sigma(1)}\right)\cdot ... \cdot \#\left(\AbBeta_{g-1} \cap
h_{g-1}\AbAlpha_{\sigma(g-1)}\right)[h_1]\cdot...\cdot[h_{g-1}] \\
&=& 
(-1)^{\epsilon}\sum_{\sigma\in\Perm{g-1}} (-1)^{\sigma} \left(\sum_{h_1\in
\TorsFree}\#\left(\AbBeta_1\cap h_1 \AbAlpha_{\sigma(1)}\right)[h_1]\right)\cdot ... \cdot
\left (\sum_{h_{g-1}\in\TorsFree}
\#\left(\AbBeta_{g-1}\cap h_{g-1} \AbAlpha_{\sigma(g-1)}\right)[h_{g-1}]\right) \\
&=&(-1)^\epsilon\tDelta_{g,g}.
\end{eqnarray*}
In the above computation, $\Perm{g-1}$ denotes the permutation group on
$g-1$ letters, and for each $\sigma\in\Perm{g-1}$, $(-1)^{\sigma}$
denotes the sign of the permutation, and
$\epsilon=\frac{g(g-1)}{2}$. 
Note that the sign comes about in
the formula for intersection number in the symmetric product: the
intersection number of $\beta_1\times...\times\beta_{g-1}$ and
$\alpha_1\times...\times\alpha_{g-1}$
 in the symmetric product is given
by $(-1)^{\frac{(g-1)(g-2)}{2}}$ times
the determinant of the matrix 
$\left(\#\left(\beta_{i}\cap
\alpha_{j}\right)\right)_{i,j}$, where $i,j=1,...,g-1$. 
\qed
\vspace{0.2in}

\subsection{Refinements}
\label{subsec:Refinements}

We discuss two refinements in the above discussion: signs, and
torsion in $H_1(Y;\Z)$. 

If we had used the maximal Abelian cover of $Y$, which corresponds to
the subgroup $$\pi_1(Y)\longrightarrow H_1(Y)=H,$$ we would have
obtained the polynomials $C_{i,j}$ used in the discussion of
Section~\ref{sec:CalcBig}. The reduction modulo torsion of these
polynomials gives the polynomials ${\underline C}_{i,j}$ used for the
Alexander polynomial. The proof of Proposition~\ref{prop:Alex}, with
the underlines removed, gives the following refinement:

\begin{prop} 
\label{prop:TuraevDelta}
The polynomial $C_{i,j}$ is obtained from the $H$-equivariant boundary
map of the maximal Abelian cover by taking the determinant
$\Delta_{i,j}$ of the $i\times j$ minor of the boundary
map $$\partial\colon (\Z[H])^g\longrightarrow (\Z[H])^g.$$
\end{prop}

This refinement is of interest, as the minor $\Delta_{i,j}$ appears in the Turaev's refinement~\cite{Turaev} of
Milnor torsion~\cite{MilnorTorsion}.  Turaev defines torsion invariant
which, for three-manifolds with $b_1(Y)>1$, takes the form of a
function $\tau_Y\in \Z[\SpinC(Y)]$. Indeed, he shows that the torsion
satisfies a formula:
\begin{equation}
\label{eq:TuraevEq}
\tau_Y\cm(1-\mu_i)(1-\nu_j)=(-1)^{g+i+j+1}\epsilon \Delta_{i,j}[\spinc]
\end{equation}
(see Equation~(4.1.a) of~\cite{Turaev}) for some apropriate sign
$\epsilon=\pm 1$ and a carefully chosen $\SpinC$ structure over
$Y$. Note that we have departed slightly from Turaev's notation: he
defines (for manifolds with $b_1>1$) an element $\tau(Y,\spinc)\in
\Z[H]$ which depends on a choice of what he calls an Euler structure,
which he shows to be equivalent to a $\SpinC$ structure. Then, the
element $\tau_Y\in\Z[\SpinC(Y)]$ defined by $$\tau_Y =
\tau(Y,\spinc)[\overline\spinc]$$ is indepenedent of the $\SpinC$
structure used in its definition. Equivalently, $\tau_Y$ is given by
$$\tau_Y=\sum_{\spinc\in\SpinC(Y)}T(\spinc)[\spinc],$$ where $T$ is
Turaev's Torsion function from \S~5
of~\cite{Turaev}. Theorem~\ref{thm:CalcBigT} is obtained easily by
comparing Equation~\eqref{eq:TuraevEq} with
Equation~\eqref{eq:CalculateTheta}.  However, to make the signs
explicit we must make explicit the signs which go into the definition
of $\theta$, then those which go into the relationship between it and
the intersection number of the tori from
Proposition~\ref{prop:Intersection}.

We now turn to the signs in the identification between $\theta$ and
the intersection numbers of Proposition~\ref{prop:Intersection} (and
Proposition~\ref{prop:IntersectionOne}). The cobordisms between the
moduli spaces $M_{u_1,u_2}$ can be thought of as cobordisms arising
from homotopies of $\ParTheta$,
$$\ParTheta_{u_1,u_2}(D,s,t)=(\Theta_{\met_s}(D)+u_1\psi(s)\alpha_i^*,
\Theta_{\met_t}(D)-u_2\psi(1-t)\beta_j^*).$$
As such, the fibers are seen to be cobordant to the fibers of the map
(for large $u_i$) 
$$\ParTheta_\infty\colon \Sym^{g-1}(\Sigma)\times [0,1]\times[0,1]
\longrightarrow S^1_{\alpha_1}\times ... \times S^1_{\alpha_g}\times
S^1_{\beta_1}\times ...\times
S^1_{\beta_g}$$
$$(D,s,t)\mapsto (\Theta_{\met_0}(D)+u_1 s \alpha_i^*,
\Theta_{\met_1}(D)-u_2 (1-t) \beta_j^*).$$
In turn, these fibers are oriented in the same manner as the fibers of
the map
$$\Sym^{g-1}(\Sigma)\longrightarrow S^1_{\alpha_1}\times ... \times S^1_{\alpha_g}\times
S^1_{\beta_1}\times ...\times
S^1_{\beta_g}/S^1_{\alpha_i}\times S^1_{\beta_j}$$
given by
$$D\mapsto
(\Theta_{\met_0}(D),\Theta_{\met_1}(D))/S^1_{\alpha_i}\times S^1_{\beta_j}.$$
The map from 
$$\left(S^1_{\alpha_1}\times ...\times \widehat {S^1_{\alpha_i}} \times ...\times S^1_{\alpha_g}\right) \times 
\left(S^1_{\beta_1}\times ...\times \widehat {S^1_{\beta_i}} \times ...\times S^1_{\beta_g}\right)$$
to the quotient torus has degree $(-1)^{i+j+g+1}$.  The preimage of
composing $\Theta$ with the map to
$S^1_{\alpha_1}\times ...\times \widehat {S^1_{\alpha_i}}
\times ...\times S^1_{\alpha_g}$, obtained by evaluating respective
holonomies, is orientation-preserving equivalent to the torus
$\Ta{i}$.
As a consequence of the above discussion, we obtain the following
precise form for the calculation of $\theta$:
\begin{equation}
\label{eq:CompThetaSigns}
\theta (1-\mu_i)(1-\nu_j)= (-1)^{i+j+g+1}{\widetilde C}_{i,j},
\end{equation}
where the tori used for ${\widetilde C}_{i,j}$ are oriented as they
are written (with respect to some consistent ordering for the
attaching circles). 

\noindent{\bf{Proof of Theorem~\ref{thm:CalcBigT}.}} 
From Proposition~\ref{prop:CalcTheta}, Equation~\eqref{eq:TuraevEq},
and Proposition~\ref{prop:TuraevDelta}, it follows that $$ \theta=
x\cm \tau_Y, $$ for some class $x\in H^2(Y;\Z)$.  Moreover, since both
$\theta$ and $\tau_Y$ are invariant under conjugation (see
Proposition~\ref{prop:StandardProperties} for $\theta$ and \S 5
of~\cite{Turaev} for $\tau_Y$), it follows that $2x=0$.

To compare signs, note that Turaev uses a slightly different sign
conventions. For example, for a chain complex whose $C_1$ is has an
oriented basis $\{a_1,...,a_g\}$ and $C_2$ is oriented by
$\{b_1,...,b_g\}$ with $\delta b_i=0$ for $i=1,...,h$, and $\delta
b_j=a_j$ otherwise, using the homology orientation induced by
$\{\alpha_1,...,\alpha_h\}$ and $\{\beta_1,...,\beta_h\}$, the sign of
the torsion over $\R$ is $(-1)^{(g-h)h}$, so Turaev's sign-refined
torsion has sign $(-1)^{(g-h)h+N(C)}=(-1)^{1+g}$ (here, $(-1)^N(C)$ is
defined in~\cite{Turaev}). On the other hand, this orientation of
$C_1\oplus C_2$ differs from the orientation induced from our
conventions by a sign of
$(-1)^{\frac{(g-1)(g-2)}{2}+\frac{(h-1)(h-2)}{2}}$. Comparing with the
sign difference between $\Delta_{i,j}$ and ${\widetilde C}_{i,j}$
it follows that 
$$\theta=(-1)^{\frac{(h-1)(h-2)}{2}} x \tau_Y.$$
\qed

\noindent{\bf{Proof of Theorem~\ref{thm:CalcOneT}.}} 
Similarly to the above, we have that $${\widetilde
C}_{i,j}=(-1)^\epsilon x \cm\Delta_{i,j}$$ The formulas relating
$\theta$ with ${\widetilde C}_{i,j}$ (see
Equation~\eqref{eq:CalculateTheta}) and $T_t$, $T_{t^{-1}}$ with
$\Delta_{i,j}$ (see \S~4 and 5 of~\cite{Turaev}) imply that
$T_t(\spinc)=\theta(\spinc+x)$ if ${\underline\spinc}$,
${\underline\spinc}+{\underline x}\geq 0$, and
$T_{t^{-1}}(\spinc)=\theta(\spinc+x)$ if ${\underline\spinc}$,
${\underline\spinc}+{\underline x}\leq 0$.  Moreover, by
Theorem~\ref{thm:CalculateZ} and the corresponding relation between
$T_t$ with Milnor torsion from~\cite{Turaev}, it follows that
${\underline\tau}'={\underline\theta}$, where ${\underline\tau}'$ is
induced from $\tau'$ in the usual manner. Since these are non-zero
polynomials, it follows that ${\underline x}=0$.  Now (in view of the
discussion of signs given in the previous proof), this implies that
$\tau'=x \theta$. Finally, since both $\tau'$ and $\theta$ are
symmetric under conjugation, it follows that $2x=0$.
\qed
\section{Wall-Crossing for $\theta$ when $b_1(Y)=1$}
\label{sec:WallCross}

The present section is meant as a technical appendix, where we show
that the definition of $\theta(\spinc)$ is independent of the
perturbation used in its definition, in the case where
$b_1(Y)=1$. Recall that this was already established in
Section~\ref{sec:DefTheta} for the case where $H_1(Y;\Z)\cong \Z$; this is
Proposition~\ref{prop:IndepOfMetOneZ}. Indeed, the arguments from
that section show that there are at most two values which
$\theta_{\eta_0\times\eta_1}(\spinc)$ can assume (for generic, small
$\eta_0\times\eta_1$), depending on the component in $H^2(Y;\R)-0$ in
which $\delta(\eta_0-\eta_1)$ lies. Thus, if we fix an
identification $H^2(Y;\R)\cong\R$, there are {\em a priori} two
invariants $\theta^\pm(\spinc)$, corresponding to the sign of
$\delta(\eta_0-\eta_1)$ under the identification.

Our goal, then, is to prove the following restatement of
Proposition~\ref{prop:IndepOfMetOne}:

\begin{prop}
\label{prop:IndepOfMetOneGen}
When $b_1(Y)=1$, then the two invariants $\theta^+(\spinc)$ and
$\theta^-(\spinc)$ agree. 
\end{prop}

In essence, this proposition amounts to the calculation of a
``wall-crossing formula'' much like the sorts of formulae one runs
across in gauge theory (see~\cite{HCobordFalse}). In the case at hand,
we have that the wall-crossing formula is trivial, which is what one
expects from the analogy with Seiberg-Witten theory, as the
perturbation is ``small'' (see~\cite{MengTaubes} for a discussion of
the three-dimensional Seiberg-Witten invariant).

To prove Proposition~\ref{prop:IndepOfMetOneGen}, we explicitly
identify the difference, in the following lemma.

\begin{lemma}
\label{lemma:Difference}
The difference $\theta^+(\spinc)-\theta^-(\spinc)$ is given by the 
intersection number
$$\theta^+(\spinc)-\theta^-(\spinc)=\#\{(t,D)\in [0,1]\times \OurSym
\big|
\OurTheta_{\met_t}(D)\in \Lag_0(\spinc)\cap \Lag_1(\spinc)\}.$$
\end{lemma}

Strictly speaking, to make sense of this intersection, we must choose
a ``generic'' allowable path of metrics $\met_t$, i.e. a path of
metrics $\met_t$ with the property that $\met_0$ and $\met_1$ are
allowable for $U_0$ and $U_1$ as usual, for which the map $[0,1]\times
\Sym^{g-1}(\Sigma)\longrightarrow \Jac$ given by $(t,D)\mapsto
\Theta_{\met_t}(D)$ is transversal to the one-manifold $\Lag_0\cap
\Lag_1\subset \Jac$. We can find such a family, according to the
following transversality result, whose proof is given in~\cite{ThetaCasson}:

\begin{theorem}
\label{thm:GenericMetrics}
If $\Sigma$ is an oriented two-manifold with genus greater than $1$,
then the $g-1$-fold Abel-Jacobi map
$$\Theta\colon \Met(\Sigma)\times \Sym^{g-1}(\Sigma)\rightarrow
\Jac$$
is a submersion.
(Here, $\Met(\Sigma)$ denotes the space of all metrics on $\Sigma$.)
\end{theorem}

In particular, standard transversality theory allows us to conclude:

\begin{cor}
Any smooth path of metrics $\met_t$ can be approximated arbitrarily
well (in $C^0$) by smooth paths 
$\met_t'$ for which the map
$[0,1]\times \Sym^{g-1}(\Sigma)\longrightarrow \Jac$
given by
$(t,D)\mapsto \Theta_{\met_t'}(D)$
is transverse to $\Lag_0\cap\Lag_1\subset \Jac$.
\end{cor}

(Note that the hypothesis that $g>1$ is not needed in the corollary;
for if $g=0,1$, then $\Lag_0\cap\Lag_1$ is automatically disjoint from
the image of $\Theta$ for any metric and, in fact,
Proposition~\ref{prop:IndepOfMetOneGen} is clear.)

The formulation given in Lemma~\ref{lemma:Difference} is useful, since
we can give the intersection number appearing there an interpretation
in terms of the index theory, from which it can be explicitly
computed. To this end, we find it convenient to use the notion of {\em
spectral flow} introduced in~\cite{APSIII}: given a one-parameter
family of self-adjoint, Fredholm operators, the spectral flow is the
intersection number of the (real) spectra of the operators with the
zero eigenvalue. We will be interested in the case where the operators
are Dirac operators coupled to $\SpinC$ connections with traceless
curvature.  Specifically, the set $\Lag_0(\spinc)\cap
\Lag_1(\spinc)$ is canonically identified with the space of gauge
equivalence classes of such connections in the $\SpinC$ structure
$\spinc$: it is empty unless $\spinc$ is torsion, in which case it can
also be identified with the circle $S^1=H^1(Y;\R)/H^1(Y;\Z)$.  (A $\SpinC$
connection is a connection on the spinor bundle $W$ of the $\SpinC$
structure and which is compatible with the Levi-Civita connection on
the tangent bundle; and the gauge group is the space of circle-valued
functions over $Y$.) The crux of the argument, then, is the following:

\begin{prop}
\label{prop:SFInterp}
The real spectral flow for the $\SpinC$ Dirac operator around the
circle $H^1(Y;\R)/H^1(Y;\Z)$, thought of as parameterizing equivalence classes
of traceless connections $A_t$ in some torsion $\SpinC$ structure
$\spinc$, is also calculated by the intersection number (with a factor
of two): $$\SF_{S^1}(Y,{A_t})=\pm 2\#\{(t,D)\in [0,1]\times
\OurSym
\big|
\OurTheta_{\met_t}(D)\in \Lag_0(\spinc)\cap \Lag_1(\spinc)\}.$$
\end{prop}

\begin{remark}
The factor of $2$ is an artifact of the complex linearity of the
$\SpinC$ Dirac operator. Moreover, the sign depends on orientation
conventions used.
\end{remark}

Proposition~\ref{prop:IndepOfMetOneGen} is an immediate consequence of
this spectral flow interpretation, together with the Atiyah-Singer
index theorem.

\vskip0.3cm
\noindent{\bf Proof of Proposition~\ref{prop:IndepOfMetOneGen}.}
A circle $[A_t]$ of gauge equivalence classes of $\SpinC$ connections
in the $\SpinC$ structure $\spinc$ over $Y$ naturally induces a
$\SpinC$ structure $\spincX$ on $X=S^1\times Y$, endowed with a (gauge
equivalence class of) $\SpinC$ connection whose restriction to the
slice $e^{it}\times Y$ is identified with $[A_t]$.  According to
Atiyah-Patodi-Singer (see~\cite{APSIII}), the spectral flow of the
Dirac operator around the circle of operators $[A_t]$ is the (real)
index of the Dirac operator on $S^1\times Y$, in the $\SpinC$
structure $\spincX$, which, according to the Atiyah-Singer index
theorem, is in turn calculated by $$ \ind \Dirac(S^1\times Y,\spincX)
=\frac{c_1(\spincX)^2}{4}-\frac{\sigma}{4}, $$ where $\sigma$ is the
signature of the intersection form of $S^1\times Y$.  In fact, the
index vanishes, since the signature $\sigma$ of $S^1\times Y$
vanishes, and the square of $c_1(\spincX)$ is also easily seen to vanish,
too, since for any fixed point $p\in S^1$, the restriction of the
$c_1(\spincX)$ to the slice $\{p\}\times Y$ is $c_1(\spinc)$, which is
a torsion class. Thus, in light of Lemma~\ref{lemma:Difference} and
Proposition~\ref{prop:SFInterp} the difference in the invariants must
vanish.
\qed
\vspace{0.2in}

We dispense first with the proof of Lemma~\ref{lemma:Difference}, and
then return to the more involved Proposition~\ref{prop:SFInterp}.

\vskip0.3cm
\noindent{\bf{Proof of Lemma~\ref{lemma:Difference}}}
Fix a path of perturbations $$\wp\colon[-1,1]\longrightarrow \NdZero \subset
\Quot_0\times\Quot_1$$ $\wp(t)=(\eta_0(t)\times\eta_1(t))$ for which
$\delta(\eta_0(t)-\eta_1(t))$ is a monotone increasing function of
$t$, which crosses $0$ at $t=0$ (here, $\NdZero$ is the neighborhood
defined in Proposition~\ref{prop:IndepOfMetBig}). According to the
transversality result (Theorem~\ref{thm:GenericMetrics}), we can find
a one-parameter family of metrics $\met_t$ so that $$\ParTheta\colon
[0,1]\times[0,1]\times
\Sym^{g-1}(\Sigma)
\longrightarrow \Quot_0\times\Quot_1$$
is transverse to $\wp_t$. 
The set 
$$\ParTheta^{-1}(\wp[-1,1])\cap \{(s,t,D)\big| s\leq t\}$$
is a one-dimensional manifold-with-boundary, whose boundary is
\begin{eqnarray*}
\#\partial \ParTheta^{-1}(\wp[-1,1])
&=& \#\ParTheta^{-1}(\wp(+1))
\cap \{(s,t,D)\big| s\leq t\} \\
&& -\#\ParTheta^{-1}(\wp(-1))
\cap \{(s,t,D)\big| s\leq t\} \\
&& - \#\ParTheta^{-1}(\wp[-1,1])\cap \{(s,s,D)\}.
\end{eqnarray*}
The points in these sets are partitioned naturally into $\SpinC$ structures.
For a fixed $\SpinC$ structure $\spinc$,
the  signed number of points 
in the first two sets calculates
$\theta^+(\spinc)$ and $\theta^-(\spinc)$ respectively while the
intersections in the third set all occur at $\wp(0)$, and indeed they
correspond to 
$$\#\{(t,D)\in [0,1]\times \OurSym
\big|
\OurTheta_{\met_t}(D)\in \Lag_0(\spinc)\cap \Lag_1(\spinc)\}.$$
The lemma follows.
\qed
\vskip0.3cm

The proof of Proposition~\ref{prop:SFInterp} uses splitting techniques
for spectral flow (see~\cite{Yoshida} and \cite{KK}): the spectral
flow around the circle $H^1(Y;\R)/H^1(Y;\Z)$ has a contribution from the
handlebodies and from the cylinder $\Sigma\times \R$. Formal
properties (reminiscent of the special case of
Proposition~\ref{prop:IndepOfMetOne} proved in
Section~\ref{sec:Introduction}) show that the contribution from the
handlebodies vanishes. The spectral flow on the cylinder is then
identified with the intersection number, in a manner akin to Yoshida's
algorithm for calculating the instanton Floer
grading~\cite{YoshidaII}.

We turn our attention, then, to the splitting of spectral flow.  We
will consider the spectral flow of the $\SpinC$ Dirac operator on
various three-manifolds $Z$, fixing the metric, and varying the
$\SpinC$ connection $A$, but keeping its curvature form to be
traceless.  For a fixed metric and $\SpinC$ structure $\spinc$, the
set of gauge equivalence classes of such connections, is analogous to
the Jacobian of a Riemann surface: it is (non-canonically) identified
with the torus $H^1(Z;S^1)$. If $\spinc$ is actually induced from a
spin structure, then this spin structure gives a canonical
identification between the two sets.

Over the cylinder $\R\times\Sigma$ given a product metric, a $\SpinC$
structure amounts to a $\SpinC$ structure on $\Sigma$, which in turn
corresponds to a line bundle $E$ over $\Sigma$ (by tensoring $E$ with
the canonical $\SpinC$ structure on $\Sigma$). Moreover, a $\SpinC$
connection corresponds to a connection on the line bundle $\R\times E$
over $\R\times\Sigma$.  With respect to the canonical splitting of the
spinor bundle over the cylinder $W=E\otimes \left(\C \oplus
\KSig^{-1}\right)$, the Dirac operator on the cylinder can be written
as $$\Dirac_A = \DDt +
\sqrt{2}\left(\begin{array}{cc} 0 & \DBar_A \\ -\DBar_A & 0
\end{array}
\right),$$
where all derivatives here mean covariant derivatives coupled to $A$
(so that the operator in the second term of the above decomposition is
the $\SpinC$ Dirac operator on $\Sigma$).
The curvature of the determinant line bundle vanishes iff the
corresponding connections on $E$ over $\Sigma$ have normalized
curvature form, in the sense of Section~\ref{sec:DefTheta}. 

Suppose that $Y$ is a three-manifold with a Heegaard splitting, which
we write as
\begin{equation}
\label{eq:MetricHeegaardDecomposition}
Y=U_0\cup_{\Sigma_0}\cup
\left([0,1]\times\Sigma\right)\cup_{\Sigma_1} U_1,
\end{equation}
where, of course, the surfaces $\Sigma_0$ and $\Sigma_1$ are topologically
identified with $\Sigma$.  A path  $\met_t$  of metrics over $\Sigma$, which is
constant near $t=0$ and $t=1$, gives rise to a metric on the cylinder
$[0,1]\times
\Sigma$, given by the formula $dt^2+\met_t$, which is product-like
near the boundary.  Fix any metric over $U_0$ (resp. $U_1$) with
boundary isometric to $\Sigma$ with metric $\met_0$ (resp. $\met_1$).
Then, these data naturally glue together to give a metric on
$Y$. 

Suppose $A$ is a $\SpinC$ connection over $Y$ with traceless
curvature, and for both $i=0,1$, the metric $\met_i$ is $U_i$
allowable. Then, the (two-dimensional) Dirac operator on the
boundaries of the three pieces of the decomposition of
$Y$ of Equation~\eqref{eq:MetricHeegaardDecomposition} have no kernel. In this
general situation, Atiyah-Patodi-Singer (see~\cite{APSI}) show that
the restriction of the Dirac operator to the three individual pieces
with APS boundary conditions is a Fredholm operator. Thus, if we have
a one-parameter family of connections $A_t$ on $Y$ whose curvature has
vanishing trace, it makes sense to speak of the spectral flow of the
Dirac operators restricted to these three pieces. Indeed, a fairly
elementary version of the splitting technology for spectral flow gives
the following:

\begin{prop}
\label{prop:SplitSF}
Let $Y$ be a three-manifold decomposed (metrically) as in
Equation~\eqref{eq:MetricHeegaardDecomposition}, with $U_i$-allowable
metrics on $\Sigma_0$ and $\Sigma_1$.  Let $[A_t]$ be any closed path
of gauge equivalence classes of $\SpinC$ connections over $Y$ with
traceless curvature.  Then, the spectral flow of the Dirac operator
coupled to the $[A_t]$ splits as a sum of the spectral flows of the
Dirac operator restricted to the three pieces (with APS boundary
conditions): $$
\SF(Y,[A_t])=\SF(U_0,[{A_t}|_{U_0}]) +
\SF([0,1]\times\Sigma ,[{A_t}|_{[0,1]\times\Sigma}]) + \SF(U_1,[{A_t}|_{U_1}]).$$
\end{prop}

The above result is standard (see for example~\cite{Bunke},
or~\cite{CLMII} for a more general result). It is proved by showing
that for metrics on $Y$ with sufficiently long cylinders $[-T,T]\times
\Sigma$ inserted around $\Sigma_0$ and $\Sigma_1$, the small
eigenmodes on $Y$ are approximated by the small eigenmodes of the
operators restricted to the individual pieces, under a splicing
construction. Since the spectral flow around the circle is independent
of these ``neck-length'' parameters (by the homotopy invariance of
spectral flow), we do not need to include them in the above statement
of the propositioon.

In fact, the only term which contributes in the above decomposition of
the spectral flow is the middle term (the one over the cylinder
$[0,1]\times\Sigma$), according to the following result, which is a
formal consequence of the conjugation action:

\begin{lemma}
\label{lemma:SFHandlebody}
Let $U$ be a handlebody which bounds the surface $\Sigma$, equipped
with a metric which is product-like near its boundary, where it
induces a $U$-allowable metric. The spectral flow of the Dirac
operator vanishes around any closed path $[A_t]$ of gauge equivalence
classes of $\SpinC$ connections, all of whose curvature is traceless.
\end{lemma}

\begin{proof}
Since $H^2(U;\Z)=0$, there is a unique $\SpinC$ structure over $U$.
Moreover, there is a complex-antilinear involution $$j\colon
W\longrightarrow W$$ of the spinor bundle which commutes with Clifford
multiplication (actually, this involution exists in much more general
contexts, and can be thought of as the basis for the conjugation
action on the set of $\SpinC$ structures described in
Section~\ref{sec:Introduction}).  It follows that if $B$ is the
connection on $W$ coming from a spin structure $\spin_0$ on $Y$,
$a\in\Omega^1(U)$, then $$\Dirac_{B+ia}(j\Psi)=j\Dirac_{B-ia}(\Psi).$$
We can express any given closed path $[A_t]$ as $[B+ia_t]$, where
$\{a_t\}$ is a one-parameter family of closed one-forms 
which induces a closed path $[a_t]$ in
$H^1(U;S^1)$; and the homotopy invariance of the spectral flow ensures
that the spectral flow of the Dirac operator around $[A_t]$ depends
only on the free homotopy class of $[a_t]\subset H^1(Y;\R)/H^1(Y;\Z)$ (in
particular, it is independent of the spin structure). Now, the
conjugation symmetry gives us that
$$\SF(U,[B+ia_t])=\SF(U,[B-ia_t]),$$ but these two spectral flows have
opposite signs: the path $[-a_t]$ is homotopic to the path $[a_t]$
given the opposite orientation.  Thus, the spectral flow around the
$[A_t]$ must vanish.
\end{proof}

\vskip0.3cm 
\noindent{\bf Proof of Proposition~\ref{prop:SFInterp}.}  

In view of Proposition~\ref{prop:SplitSF} and
Lemma~\ref{lemma:SFHandlebody}, the spectral flow over $Y$ is
determined by the spectral flow around a loop over the cylinder
$[0,1]\times \Sigma$. So we turn our attention now to the study of the
spectral flow over a cylinder. In fact, it is useful to consider a
more general setting -- spectral flow along a not necessarily closed
path of operators on the cylinder $[0,1]\times\Sigma$.

Let $(\met_t, A_t)_{t\in[0,1]}$ be a path of metrics over $\Sigma$ and
connections $[A_t]\in\Jac_{\met_t}$, both of which are constant near
the $t=0,1$ endpoints. We can canonically extend the paths $(h_t,A_t)$
for all $t\in\R$, so that they remain constant for $t\leq 0$ and
$t\geq 1$.  Suppose, now, that $A_0$ does not lie in the
$\met_0$-theta divisor and similarly, $A_1$ does not lie in the
$\met_1$-theta divisor. Then, the associated Dirac operator on
$[0,1]\times
\Sigma$ -- the one for the metric $dt^2+\met_t$ and the spin connection
obtained by viewing the path of connections $\{A_t\}$ as a single
connection over $[0,1]\times \Sigma$ -- is a Fredholm operator on
$L^2$.  Suppose moreover that each pair $(A_t,\met_t)$ for $t\in[0,1]$
misses the $\met_t$-theta divisor. Then, if we rescale the family in
the $\R$ direction to move sufficiently slowly, then the Dirac
operator on the cylinder $\R\times\Sigma$ will have no kernel (this
follows from a standard adiabatic limit argument, a proof is given in
Proposition~\ref{Casson:prop:StretchOutKernel} of~\cite{ThetaCasson}).
Moreover, if we have a two-parameter family: $$H\colon
[0,1]\times[0,1]\rightarrow
\Met(\Sigma)\times \Jac,$$ where the boundary misses the theta
divisor, then we get a one-parameter family of self-adjoint, Fredholm operators
$D(s)$
indexed by $s\in[0,1]$ which we get from $H(s,t)$ by fixing the $s$
coordinate and allowing $t$ to vary. According to
the adiabatic limit statement, the spectral flow vanishes if $H$
always misses the theta divisor. Indeed, by the homotopy invariance of
the spectral flow and the connectedness of the theta
divisor, the spectral flow depends only on the homological
intersection number of $H$ with the theta divisor.  This proves that
\begin{equation}
\label{eq:SFCyl}
\SF(D(s))=\mu \cm \# (H\cap\Theta),
\end{equation}
for some integer $\mu$, which {\em a priori} depends only on the
genus $g$ of $\Sigma$ (which we suppress from the notation wherever it
is convenient).

To determine $\mu$, we consider a simple model case: we construct a
two-parameter family of metrics and connections which intersects the
theta divisor once, extend it naturally over a three-manifold,
calculate the Chern class of the extension, and then compare with the
result obtained from the Atiyah-Singer index theorem to calculate the
spectral flow.  View the surface $\Sigma$ as a connected sum of $g$
tori $F_1,...,F_g$, and let $\{\alpha_1,...,\alpha_g\}$ be a complete
set of attaching circles with $\alpha_i$ supported in $F_i$. Also, for
$i=1,...g$, let $\beta_i$ be a simple closed curves in $F_i$ so that
$\#(\alpha_i\cap\beta_i)=1$.  Consider a two-parameter family
$$H\colon [0,1]\times[0,1]\longrightarrow \Met(\Sigma)\times
\Jac$$ where the metric $\met$ is held constant (we will say how it is chosen
in a moment), and the holonomy of the connection associated to
$H(s,t)$ around $\alpha_g$ is $e^{2\pi i t}$ (independent of $s$), the
holonomy of $A$ around $\beta_g$ is $e^{2\pi i s}$, and all the other
holonomies are trivial (here, the holonomies are measured relative to
a spin structure on $\Sigma$ which extends to the handlebody
determined by $\{\alpha_1,...,\alpha_g\}$).  

For metrics which are sufficiently stretched out along the connected
sum curves and the attaching curves $\{\alpha_i\}_{i=1}^g$, the image
of $H$ intersects the theta divisor transversally in a single point
(with some appropriate choice of sign): this fact is closely related
to Corollary~\ref{cor:HolonomyConstraints}. To see this, as in the
proof of that proposition, one uses Proposition~\ref{prop:RangeSplice}
to conclude that (for metrics $\met$ on $\Sigma$ which are
sufficiently stretched out) the intersection is contained in the image
of a splicing map, which is $C^1$ close to a map $$\FComp_1 \times
... \times
\FComp_{g-1} \longrightarrow H^1(F_1;S^1)\times ... \times
H^1(F_{g-1};S^1)\times H^1(F_g;S^1),$$ which is the degree one
Abel-Jacobi map on the first $g-1$ torus factors, and constant on the
final factor. Requiring the holonomies to be trivial around the
$\alpha_i$ and $\beta_i$ for $i=1,...,g-1$, is equivalent to
retricting to a single point in the domain. Since the degree one
Abel-Jacobi map is a diffeomorphism, it follows that $H$ indeed
intersects the theta divisor transversally in a single point.

For each fixed $s$, the family of connections $H(s,t)$, where $t$
varies, canonically extends as a flat connection over (at both $t=0$
and $t=1$) the handlebody $U_0$ obtained by surgeries
along the $\alpha_i$, to give connections $A_s$ on a line bundle $L$
over the three manifold $Y$ obtained as the $g$-fold connected sum
$$Y=\overbrace{\left(S^1\times S^2\right)\#...\#
\left(S^1\times S^2\right)}^g.$$ From the construction of $A_s$, it is
clear that its curvature vanishes on all but the final connected
summand. Indeed, it is easy to see that the first Chern class $c_1(L)$
is dual to the two-sphere in that summand.  Moreover, the connections
at $s=0$ and $s=1$ are gauge equivalent, via a gauge equivalence which
extends over $Y$. Letting $u$ denote the gauge transformation over
$Y$. Note that the gauge transformation is non-trivial only over the
final connected summand of $Y$, where it gives a map of degree one on
its circle $\beta_g$. Now, the $A_s$ naturally induce a connection on
the line bundle $M$ over $S^1\times Y$ obtained from $[0,1]\times L$
by identifying $\{0\}\times L$ with $\{1\}\times L$ using the gauge
transformation $u$.  From what we know about $L$ and $u$, it follows
easily that the first Chern class of the line bundle $M$ is Poincar\'e
dual to $S^1\times \beta_g$ plus the sphere $S^2$ which appears in the
final connected summand.  Thus, tensoring any spin structure over
$S^1\times Y$ with $M$, we obtain a $\SpinC$ structure $\spincX$ whose
first Chern class is twice the first Chern class of $M$, so according
to the Atiyah-Singer index theorem, the index of the Dirac operator
coupled to $L$ is given by $$\ind
\Dirac(S^1\times Y,\spincX)= 2.$$ Note that this index calculates the
spectral flow around the $S^1$-factor, which consists only of the
contribution of the cylinder (according to
Proposition~\ref{prop:SplitSF} and Lemma~\ref{lemma:SFHandlebody}).
Since the intersection number of our family with the theta divisor
consisted of a single, isolated point, it follows that $\mu=\pm 2$ in
Equation~\eqref{eq:SFCyl} (in particular, $\mu$ is independent of the
genus $g$).
\qed 
\vskip0.3cm

\commentable{
\bibliographystyle{plain}
\bibliography{biblio}
}

\end{document}